\documentclass
{amsart}
\usepackage{graphicx}
\newtheorem{thm}{Theorem}[section]
\theoremstyle{definition}
\newtheorem{cor}[thm]{Corollary}
\newtheorem{prop}[thm]{Proposition}
\newtheorem{defn}[thm]{Definition}
\newtheorem{lem}[thm]{Lemma}

\newtheorem{note}[thm]{Notation and Remark}
\newtheorem{rem}[thm]{Remark}

\newtheorem{ex}[thm]{Example}
\newtheorem{Conclusion}[thm]{Conclusion}
\newtheorem{exs}[thm]{Examples}
\newtheorem{co}[thm]{Question}
\numberwithin{equation}{section}
\begin{document}
\thispagestyle{empty}
\vspace{0.8cm}
{\begin{center}
{\large \bf Survey on second submodules of modules over commutative rings

} \vspace*{3mm}


{\bf Faranak Farshadifar\footnote{Department of Mathematics Education, Farhangian University, P.O. Box 14665-889, Tehran, Iran.\ E-mail: {f.farshadifar@cfu.ac.ir}}}
\end{center}


\begin{abstract}
Let $R$ be a commutative ring with identity.
The concept of second submodule of an $R$-module (as a dual notion of prime submodules)  was introduced and studied by  S.Yassemi in 2001. This notion has obtained a great attention by many authors and now there is a considerable amount of research concerning this class of modules. The main purpose of this paper is to collect these results and provide a useful source for those who are interested in research in this field.
\end{abstract}


\tableofcontents

\section{Introduction}
\noindent
Throughout this paper, $R$ will denote a commutative ring with
identity and ``$\subset$" will denote the strict inclusion. Further, $\Bbb Z$ and $\Bbb N$ will denote the ring of integers and the set
positive integers, respectively. We use $N\leq M$ to indicate that $N$ is a submodule of a module $M$.
For a submodule $N$ of an $R$-module $M$, the colon ideal of $M$ into $N$, $(N:_RM)=\{r \in R: rM \subseteq N\}=Ann_R(M/N)$.
For any unexplained notions or terminology please see \cite{AF74, MR2455920,NV82}.

A proper submodule $P$ of an $R$-module $M$ is said to be \textit{prime} if for any
$r \in R$ and $m \in M$ with $rm \in P$, we have $m \in P$ or $r \in (P :_R M)$ \cite{MR183747, Da78}.
In \cite{Y01}, S. Yassemi introduced the dual notion of prime submodules (i.e., second submodules) over a commutative ring and investigated the first properties of this class of modules. Annin \cite{Ann2002} called
these \emph{coprime modules} (see also \cite{Wij2006}) and used them to
dualize the notion of attached primes. A non-zero submodule $S$ of $M$ is said to be \textit{second} if for each $a \in R$, the homomorphism $ S \stackrel {a} \rightarrow S$ is either surjective or zero \cite{Y01}. The set of all second submodules of $M$ is called the \textit{second spectrum of $M$} and denoted by
$Spec^s(M)$ \cite{AF21}.

Theory of prime ideals is an important tool in classical algebraic geometry.
In development of algebraic geometry and module theory, some generalizations
for the concept of prime ideals has arisen. There are extensive
investigations concerning the prime submodules of modules and this notion
has attracted attention by a number of authors. It is natural to ask the following
question: To what extent does the dual of these results hold for second
submodules of an R-module? To answer this question, there are a number of authors who have studied the second submodules. Now there is plenty of useful information which has been investigated by many authors \cite{Y01,AF11, AF12, AF21, AF25, CAS103, FF2024,   AF1089891, MR4401391, HM2, MR4049591,  MR4401392}. The main purpose of this paper is to collect these results and provide a useful source for those who are interested in research in this field. We refer the reader to  \cite{2026F} for a survey on the second spectrum of a module.
\section{Second submodules}
Let $M$ be an $R$-module.
A submodule $N$ of $M$ is said to be \textit{cocyclic} if $N\subseteq E(R/m)$ for some maximal ideal $m$ of $R$ (here $E(R/m)$ denote the injective envelope of $R/m$) \cite{MR1654522}.
A proper submodule $N$ of $M$ is said to be \emph{completely
irreducible} if $N=\bigcap _ {i \in I}N_i$, where $ \{N_i \}_{i \in I}$ is a family of submodules of $M$, implies that
$N=N_i$ for some $i \in I$. It is easy to see that every submodule of $M$ is an intersection of completely irreducible submodules of $M$. Thus the intersection of all completely irreducible submodules of $M$ is zero  \cite{FHo06}.

\begin{lem}\label{t222.4}\cite[Remark 1.1]{FHo06}
A submodule $L$ of an $R$-module $M$ is completely irreducible if and only if $M/L$ is a cocyclic module .
\end{lem}

Let $\mathfrak{p}$ be a prime ideal of $R$ and let $N$ be a submodule of an $R$-module $M$. Then $N^{ec}= \{ m \in M: cm \in N $ for some $c  \in R\backslash \mathfrak{p} \}$ and  it is called the $\mathfrak{p}$-closure of $N$ and denoted by $cl_{\mathfrak{p}}(N)$ \cite{MS93,MS08}.

\begin{defn}\cite[Definition 2.7]{AF11})
Let $\mathfrak{p}$ be a prime ideal of $R$ and let $N$ be a
submodule of an $R$-module $M$. \emph{The $\mathfrak{p}$-interior of $N$ relative to $M$}
is defined as the set
$$
I^M_{\mathfrak{p}}(N)= \cap \{L \mid  L \\\ is \\\ a \\\ completely \\\
irreducible
\\\ submodule \\\ of \\\ M\\\ and
$$
$$
 rN\subseteq L \\\ for \\\ some \\\ r \in R-\mathfrak{p} \}.
$$
This can be regarded as a dual notion of the $\mathfrak{p}$-closure of $N$. It is clear that $\cap _{r \in R-\mathfrak{p}}rN \subseteq I^M_{\mathfrak{p}}(N) \subseteq
N$.
\end{defn}

\begin{defn}\label{d2.1}\cite[Definition 1.1]{Y01}
A non-zero submodule $S$ of an $R$-module $M$ is said to be \emph{second} if for each $a \in R$, the homomorphism $ S \stackrel {a} \rightarrow S$ is either surjective or zero. This implies that $Ann_R(S) = \mathfrak{p}$ is a prime ideal of $R$, and $S$ is said to be $\mathfrak{p}$-second.
\end{defn}

Let $M$ be an $R$-module. The dual notion of $Z_R(M)$, the
set of zero divisors of $M$ \cite{Y0551}, is denoted by $W_R(M)$ and defined by
$$
W_R(M)= \{ a \in R: aM \not =M \}.
$$

A prime ideal $\mathfrak{p}$ of $R$ is said to be a \textit{weakly coassociated prime} of an $R$-module $M$ if there
exists a cocyclic homomorphic image $L$ of $M$ such that $\mathfrak{p}$ is a minimal element in
$V(Ann_R(L)) = \{q \in Spec(R)|Ann_R(L) \subseteq q\}$. The set of weakly coassociated prime
ideals of $M$  is denoted by $\widetilde{Coass}(M)$\cite{Y0551}.

The $R$-module $M$ is said to be \textit{divisible} if for each $x \in M$ and for each non-zero
divisor $r \in R$ there exists $y \in M$ such that $ry = x$ \cite{MR2455920}.
\begin{thm}\label{d2.1}\cite[Theorem 1.3]{Y01} and \cite[Theorem 2.10]{AF12}
Let $S$ be a non-zero submodule of an $R$-module $M$, with $Ann_R(S) = \mathfrak{p}$. Then the
following are equivalent:
\begin{itemize}
\item [(a)]  $S$ is a $\mathfrak{p}$-second submodule of $M$;
\item [(b)]  $S$ is a divisible $R/{\mathfrak{p}}$-module;
\item [(c)]  $rS = S$ for all $r \in R\setminus \mathfrak{p}$;
\item [(d)]  $IS = S$ for all ideals $I \not \subseteq \mathfrak{p}$;
\item [(e)]  $W_R(S) = \mathfrak{p}$;
\item [(f)]  $\widetilde{Coass} (S) = \{\mathfrak{p}\}$;
\item [(g)] $S \not =0$ and $rS\subseteq K$, where $r \in R$ and $K$ is a submodule of $M$, implies either that $rS=0$ or $S\subseteq K$;
\item [(h)] $S \not =0$ and $rS\subseteq L$, where $r \in R$ and $L$ is a completely irreducible submodule of M, implies either that $rS=0$  or $S\subseteq K$.
\item [(i)] $\mathfrak{p}=Ann_R(S)$ is a prime ideal of $R$ and and $I^M_\mathfrak{p} (S) = S$.
\end{itemize}
\end{thm}

\begin{prop}\label{df2.1}\cite[Proposition 1.4]{Y01}
If $S$ is a submodule of an $R$-module $M$ and $Ann_R(S) = \mathfrak{m} \in Max(R)$,  then $S$
is a $\mathfrak{m}$-second submodule.
\end{prop}

\begin{prop}\label{df2.3}\cite[Proposition 1.6]{Y01}
If $S$ is a minimal submodule of an $R$-module $M$, then $S$ is a second submodule.
\end{prop}

Proposition \ref{df2.3} shows that any minimal submodule is second submodule.
So it is good to know, which modules have minimal submodules. The
next theorem gives these modules.

\begin{thm}\label{df2.5}\cite[Theorem 1.8]{Y01}
The $R$-module $M$ has minimal submodules if and only if there
exists non-zero finitely cogenerated submodule $L$ of $M$.
\end{thm}

\begin{cor}\label{df2.6}\cite[Corollary 1.9]{Y01}
If $M$ is a finitely cogenerated $R$-module, then every non-zero submodule
of $M$ contains a simple, hence second, submodule.
\end{cor}

A non-zero $R$-module $M$ is said to be \emph{secondary} if for each $a \in R$ the endomorphism of $M$ given by multiplication by $a$ is either surjective or nilpotent \cite{M1973}.
\begin{prop}\label{df2.77}\cite[Proposition 1.10]{Y01}
The following hold;
\begin{itemize}
\item [(a)] Let $S$ be a secondary submodule of an $R$-module $M$. Then $S$ is second if and
only if $Ann_R(S)\in Spec(R)$.
\item [(b)] Let $K$ be a submodule of a $\mathfrak{p}$-second module $M$. Then $K$ is a
$\mathfrak{p}$-secondary submodule if and only if $K$ is a $\mathfrak{p}$-second submodule.
\end{itemize}
\end{prop}

A submodule $N$ of an $R$-module $M$ is said to be a \textit{a minimal $\mathfrak{p}$-secondary (resp. $\mathfrak{p}$-second)}
submodule of $M$ if $N$ is a $\mathfrak{p}$-secondary (resp. $\mathfrak{p}$-second) submodule which is not
strictly contains any other $\mathfrak{p}$-secondary (resp. $\mathfrak{p}$-second) submodule of $M$\cite{Y01}.
\begin{thm}\label{df2.7}\cite[Theorem 1.11]{Y01}
The submodule $N$ of an $R$-module $M$ is minimal $\mathfrak{p}$-secondary if and only if $N$
is a minimal $\mathfrak{p}$-second submodule of $M$.
\end{thm}

\begin{defn}\label{df2.8}\cite[Definition 2.1.]{Y01}
Let $M$ be an $R$-module. Then
\begin{itemize}
\item [(a)] We say $M$ is a \textit{prime module} if the zero submodule of $M$ is a prime
submodule of $M$.
\item [(b)] We say $M$ is a \textit{second module} if $M$ is a second submodule of itself.
\end{itemize}
\end{defn}

Note that $N$ is a prime submodule of $M$ if and only if $M/N$ is a prime module.
In addition, the ring $R$ is a prime $R$-module if and only if $R$ is an integral domain.
Also, $R$ is a second $R$-module if and only if $R$ is a field \cite{Y01}.

\begin{prop}\label{df2.9}\cite[Proposition 2.2]{Y01}
Let $\mathfrak{p} \in Spec(R)$. Then the following hold:
\begin{itemize}
\item [(a)] The sum of $\mathfrak{p}$-second modules is a $\mathfrak{p}$-second module.
\item [(b)] Every product of $\mathfrak{p}$-second module is a $\mathfrak{p}$-second module.
\item [(c)] Every non-zero quotient of a $\mathfrak{p}$-second module is likewise $\mathfrak{p}$-second.
\end{itemize}
\end{prop}

\begin{thm}\label{df2.10}\cite[Theorem 2.3]{Y01}
Let $M$ be a prime module. Then the following are equivalent:
\begin{itemize}
\item [(a)] $M$ is a second module,
\item [(b)] $M$ is an injective $R/Ann_R(M)$-module.
\end{itemize}
\end{thm}

\begin{thm}\label{df2.11}\cite[Theorem 2.4]{Y01}
Let $M$ be a second module. Then the following are equivalent:
\begin{itemize}
\item [(a)] $M$ is a prime,
\item [(b)] $M$ is a flat $R/Ann(M)$-module.
\end{itemize}
\end{thm}

\begin{thm}\label{df2.12}\cite[Theorem 2.5]{Y01}
Let $M$ be a finitely generated $R$-module. If $M$ is a second module,
then $M$ is a prime module.
\end{thm}

The next result is a dual of Theorem \ref{df2.12} in a certain sense.

\begin{thm}\label{df2.13}\cite[Theorem 2.6]{Y01}
Let $M$ be an Artinian R-module. If $M$ is a prime module, then $M$
is a second module.
\end{thm}

\begin{cor}\label{df2.14}\cite[Corollary 2.7]{Y01}
Let $M$ be a finitely generated and Artinian module. Then $M$ is a
prime module if and only if $M$ is a second module.
\end{cor}

In \cite{Y0551} some functorial results for the category of modules such that
the zero submodule has a primary decomposition (or being zero) and the category
of modules that has secondary representation (or being zero) are given. The following theorem
bring the similar results for the class $\mathcal{P}$ of all prime R-modules (or being zero),
and the class $\mathcal{S} $ of all second $R$-module (or being zero).
\begin{thm}\label{df2.15}\cite[Theorem 3.1]{Y01}
Let $T$ be a linear functor over the category of $R$-modules. Then the following hold:
\begin{itemize}
\item [(a)] If $T$ is a left exact and covariant and if $M \in \mathcal{P}$, then $T(M) \in \mathcal{P}$.
In particular, if $F$ is a flat $R$-module and $M \in \mathcal{P}$, then $M\otimes F \in \mathcal{P}$, and if $M \in \mathcal{P}$, then $Hom_R(N,M) \in \mathcal{P}$ for any $R$-module $N$.
\item [(b)]  If $T$ is right exact and covariant and if $M \in \mathcal{P}$, then $T(M)\in \mathcal{S}$. In particular, if $E$ is an injective $R$-module and $M \in \mathcal{P}$, then $Hom_R(M,E) \in \mathcal{S}$.
\item [(c)] If $T$ is right exact and covariant and if $M \in \mathcal{S}$,  then $T(M) \in \mathcal{S}$. In particular, if $M \in \mathcal{S}$ then $M\otimes N \in \mathcal{S}$ for any $R$-module $N$. In addition, if $P$ is a projective $R$-module and $M \in \mathcal{S}$, then
$Hom_R(P,M) \in \mathcal{S}$.
\item [(d)] If $T$ is left exact and contravariant and if $M \in \mathcal{S}$, then $T(M) \in \mathcal{P}$. In particular, if $M \in \mathcal{S}$, then $Hom_R(M,N) \in \mathcal{P}$ for any $N$.
\end{itemize}
\end{thm}

\begin{lem} \label{l3.14}\cite[Lemma 2.13]{AF12} Let $E$ be an injective cogenerator of $R$ and
let $N$ be a submodule of an $R$-module $M$. Then $N$ is a prime
submodule of $M$ if and only if $Hom_R(M/N,E)$ is a second $R$-module.
\end{lem}

\begin{prop}\label{8l3.14}\cite[Proposition 2.12]{AF12}
Let $M$ be an $R$-module. If every non-zero submodule of $M$
is second, then for each submodule $K$ of $M$ and each ideal $I$ of $R$, we have
$(K :_M I) = (K :_M I^2)$. Also for any two ideals $A,B$ of $R$, $(K :_M A)$ and
$(K :_M B)$ are comparable.
\end{prop}

We say that a submodule $N$ of an $R$-module $M$
is \emph{coidempotent} if $N=C(N^2)$.
Also, an $R$-module $M$ is said to be \emph{fully coidempotent}
if every submodule of $M$ is coidempotent \cite[3.1]{AF122}.
\begin{thm}\label{8lfff3.14}\cite[Theorem 3.9]{AF122}
Let $M$ be a fully coidempotent $R$-module. Then every second submodule of $M$ is a minimal submodule of $M$.
\end{thm}

An $R$-module $M$ is said to be a \emph{comultiplication module} if for every submodule $N$ of $M$ there exists an ideal $I$ of $R$ such that $N=(0:_MI)$ \cite{AF07}.
An $R$-module $M$ is  said to be a \emph{weak
multiplication module} if $M$ does not have any prime submodule or for every prime submodule $N$ of $M$ there exists an ideal $I$ of $R$ such that $N=IM$ \cite{Sa95}.

\begin{defn} \label{d3.17}\cite[Definition 3.1]{AF12} We say that an $R$-module $M$ is a
\emph{weak comultiplication module} if $M$ does not have any
second submodule or for every second submodule $S$ of $M$, we have
$S=(0:_MI)$, where $I$ is an ideal of $R$.
This can be regarded as a dual notion of the weak
multiplication module.
\end{defn}
\begin{rem} \label{r3.18}\cite[Remark 3.2]{AF12} It is clear that every comultiplication
$R$-module is a weak comultiplication $R$-module. However in
general the converse is not true. For example, the $\Bbb Z$-module
$\Bbb Q$ is a weak comultiplication module because $\Bbb Q$ as a
$\Bbb Z$-module has no second submodule except $\Bbb Q$. But since $(0:_{\Bbb Q}Ann_{\Bbb Z}(\Bbb Z))= \Bbb
Q \not = \Bbb Z$, $\Bbb Q$ is not a comultiplication $\Bbb Z$-module.
\end{rem}
\begin{lem}\label{l3.19}\cite[Lemma 3.3]{AF12} Let $M$ be an $R$-module. Then we have the following.
\begin{itemize}
  \item [(a)] $M$ is a weak comultiplication module if and only if
  $S=(0:_MAnn_R(S))$ for each second submodule $S$ of $M$.
  \item [(b)] If $M$ is a weak comultiplication module, then every
  submodule of $M$ is a weak comultiplication module.
\end{itemize}
\end{lem}

For an $R$-module $M$, $Coass_R (M)$ denotes the set of all
prime ideals $\mathfrak{p}$ of $R$ such that there exists a cocyclic homomorphic image $L$ of $M$ with $Ann_R(L)=\mathfrak{p}$ \cite{Y97}.
\begin{thm} \label{t3.20}\cite[Theorem 3.4]{AF12} Let $M$ be an $R$-module. Then we have the following.
\begin{itemize}
    \item [(a)] If $M$ is a weak comultiplication $R$-module and has
    finite length, then every second submodule of
    $M$ is minimal.
    \item [(b)] If $M$ is a Noetherian weak comultiplication
    $R$-module, then $M$ has a finite number of second submodules
\end{itemize}
\end{thm}

\begin{lem} \label{l3.21}\cite[Lemma 3.5]{AF12} Let $R$ be a Noetherian ring and let $M$ be a
finitely generated $R$-module. Then we have the following.
\begin{itemize}
  \item [(a)] If $S$ is a multiplicatively closed
  subset of $R$ and $N$ is a second submodule of $M$ such that
  $Ann_R(N) \cap S= \emptyset$, then $S^{-1}N$ is a second submodule
  of $S^{-1}M$.
  \item [(b)] If for every maximal ideal $\mathfrak{p}$ of $R$, $M_\mathfrak{p}$ is a
  weak comultiplication $R_\mathfrak{p}$-module, then $M$ is a weak
  comultiplication $R$-module.
\end{itemize}
\end{lem}

\begin{thm} \label{t3.22}\cite[Theorem 3.6]{AF12} Let ($R$, $\mathfrak{p}$) be a Noetherian local ring and
let $M$ be a finite length weak comultiplication $R$-module.
Then $M$ is a comultiplication $R$-module.
\end{thm}

\begin{lem}\label{l2013.1}\cite[Lemma 1.1]{CAS13}
 Let $R$ be a commutative ring such that every prime ideal is maximal. Then the module $M$ is second if and only if $M$ is
homogeneous semisimple.
\end{lem}

\begin{cor}\label{l2013.2}\cite[Corollary 1.2]{CAS13}
Let $R$ be a commutative von Neumann regular ring. Then a nonzero
$R$-module $M$ is a second module if and only if $M$ is homogeneous semisimple.
\end{cor}

\begin{prop}\cite[Proposition 2.1]{MR3307390}
Let $M$ be a module over a commutative Noethenan ring $R$
and $S$ be a multiplicatively closed subset of $R$. If all elements of $S$ act bijectively on
$M$, then $M$ is a second $R$-module if and only if $S^{-1}M$ is a second $S^{-1}R$-module.
\end{prop}

\begin{prop}\label{pp.22}\cite[Proposition 2.2]{MR3307390}
Let $R$ be a Dedekind domain. Then every second $R$-module
is either homogeneouss emisimple or divisible.
\end{prop}

The Example \cite[Example 2.3]{MR3307390} shows that the condition on the ring $R$ in Proposition \ref{pp.22} is
necessary. This example also shows hat the converse of \cite[Corollary 2.6]{CAS13} is not
true in general.
\section{Maximal second submodules and $\mathfrak{p}$-interior}
\begin{defn}\label{d3.2}\cite[Definition 2.1.]{AF11}
 We say that a second submodule $N$ of an $R$-module $M$
is a \emph {maximal second submodule} of a submodule
$K$ of $M$, if $N \subseteq K$ and there does not exist a second
submodule $L$ of $M$ such that $N \subset L \subset K$.
\end{defn}

The following lemma can be proved easily by using Zorn's Lemma.

\begin{lem}\label{l3.3}\cite[Lemma 2.2]{AF11}
Let $M$ be an $R$-module. Then every second
submodule of $M$ is contained in a maximal second submodule of $M$.
\end{lem}

Let $M$ be an $R$-module.  A family $\{N_i \}_{i \in I}$ of submodules of $M$ is said to be an \emph{inverse family of submodules of $M$} if the intersection of two of its submodules again contains a module in $\{N_i \}_{i \in I}$. Also,
$M$ \emph {satisfies the Grothendieck's condition $AB5^*$} (\emph {the property $AB5^*$} in short) if for every submodule $K$ of $M$ and every inverse family $\{N_i\}_{i \in I}$ of submodules of $M$, $K+\cap _{i \in I}N_i= \cap _{i \in I}(K+N_i)$. Artinian and uniserial modules are examples of modules which satisfies the property $AB5^*$ \cite[p.435]{W91}.

\begin{thm}\label{t3.4}\cite[Theorem 2.3]{AF11}
 Let $M$ be a finitely cogenerated
comultiplication $R$-module which satisfies the property $AB5^*$.
Suppose that for each maximal second submodule $K$ of $M$, we have
$M/K$ is Artinian. Then the number of maximal
second submodules of $M$ is finite.
\end{thm}

\begin{thm}\label{t3.11}\cite[Theorem 2.8]{AF11}
Let $N$ be a submodule of an $R$-module $M$ such that $Ann_R(N)=\mathfrak{p}$ is a prime ideal of $R$. If
$M/I_\mathfrak{p}(N)$ is a finitely cogenerated $R$-module, then
$I_\mathfrak{p}(N)$ is a maximal $\mathfrak{p}$-second submodule of $N$.
\end{thm}

\begin{lem}\label{l3.13}\cite[Lemma 2.9]{AF11}
Let $\mathfrak{p}$ be a prime ideal of $R$ and let $M$
be an $R$-module such that $M/I_\mathfrak{p}((0:_M\mathfrak{p}))$ is a finitely
cogenerated $R$-module. If $I_\mathfrak{p}((0:_M\mathfrak{p})) \not =0$, then
$I_\mathfrak{p}((0:_M\mathfrak{p}))$ is a maximal $\mathfrak{p}$-second submodule of $(0:_M\mathfrak{p})$.
\end{lem}

\begin{cor} \label{c3.15}\cite[Corollary 2.10]{AF11} Let $\mathfrak{p}$ be a prime ideal of $R$ and let $M$ be
an $R$-module such that $M/I_\mathfrak{p}((0:_M\mathfrak{p}))$ is a finitely cogenerated
$R$-module. Then the following statements are equivalent.
\begin{itemize}
  \item [(a)] $Ann_R((0:_M\mathfrak{p}))=\mathfrak{p}$.
  \item [(b)] $I_\mathfrak{p}((0:_M\mathfrak{p}))$ is a second submodule of $M$.
  \item [(c)] There exists a second submodule $K$ of $M$ such that
  $\mathfrak{p}=Ann_R(K)$.
  \item [(d)] $I_\mathfrak{p}((0:_M\mathfrak{p})) \not = 0$.
\end{itemize}
\end{cor}

\begin{thm}\label{l62.3}\cite[Theorem 2.4]{MR3588217}
Let $\mathfrak{p} \in Spec(R)$ and $N$ be a submodule of an $R$-module $M$. Then we have the following.
\begin{itemize}
  \item [(a)] If $M$ is an Artinian $R$-module, then $I^M_\mathfrak{p}(I^M_\mathfrak{p}(N))=I^M_\mathfrak{p}(N)$.
  \item [(b)] If $M$ is an Artinian $R$-module, then $Hom_R(R_\mathfrak{p},I^M_\mathfrak{p}(N))=Hom_R(R_\mathfrak{p},N)$.
  \item [(c)] $Ann_R(N)\subseteq cl_\mathfrak{p}(Ann_R(N))\subseteq Ann_R(I^M_\mathfrak{p}(N))$.
  \item [(d)] If $M$ is an Artinian $R$-module, then $Ann_R(I^M_\mathfrak{p}(N))=cl_\mathfrak{p}(Ann_R(I^M_\mathfrak{p}(N))$.
\end{itemize}
\end{thm}

\begin{defn}\label{d2.1}\cite[Definition 2.5]{MR3588217}
We say that a submodule $N$ of an $R$-module $M$ is \emph{cotorsion-free with respect to (w.r.t.)} $\mathfrak{p}$ if $I^M_\mathfrak{p}(N)=N$, where $\mathfrak{p} \in Spec(R)$.
\end{defn}

\begin{lem}\cite[Lemma 2.6]{MR3588217}
Let $N$ ba a submodule of an $R$-module $M$ and $\mathfrak{p} \in Spec(R)$. If $N$ is
cotorsion-free w.r.t. $\mathfrak{p}$, then $N$ is cotorsion-free w.r.t. $Q$ for each $Q \in V (P)$.
\end{lem}

\begin{ex}\cite[Example 2.7]{MR3588217}
\begin{itemize}
\item [(1)] If $\mathfrak{p}\in Spec(R)$, then every $\mathfrak{p}$-secondary submodule of an $R$-module $M$ is cotorsion-free w.r.t. $\mathfrak{p}$ by \cite[Theorem 2.8]{AF124}.
\item [(2)] The $\Bbb Z$-module $\Bbb Z_{p^{\infty}}$ is cotorsion-free w.r.t. $(0)$.
\end{itemize}
\end{ex}

\begin{cor}\cite[Corollary 2.8]{MR3588217}
Let $\mathfrak{p} \in Spec(R)$ and $N$ be a submodule of an $R$-module $M$. If $N$ is
cotorsion-free w.r.t. $\mathfrak{p}$, then $Ann_R(I^M_P(N)) = S_P (Ann_R(I^M_P (N))).$
\end{cor}

\begin{thm}\label{t2.4}\cite[Theorem 2.9]{MR3588217}
Let $\mathfrak{p} \in Spec(R)$ and $N$ be a submodule of an Artinian $R$-module $M$. Then we have the following.
\begin{itemize}
  \item [(1)] $Ann_{R_P}(Hom_R(R_P,N))=(Ann_R(I^M_P(N)))_P$.
  \item [(2)] The following statements are equivalent.
  \begin{itemize}
    \item [(a)] $Hom_R(R_P,N) \not =0$.
    \item [(b)] $Ann_R(I^M_P(N))\subseteq P$.
    \item [(c)] $I^M_P(N) \not =0$.
    \item [(d)] $\mathfrak{p} \in Cosupp_R(N)$.
  \end{itemize}
\end{itemize}
\end{thm}

\begin{thm}\label{d2.1}\label{t2.6}\cite[Theorem 2.11 and Corollary 2.17]{MR3588217}
Let $\mathfrak{p} \in Spec(R)$ and $0\not =N$ be a submodule of an Artinian $R$-module $M$. Then the following statements are equivalent.
\begin{itemize}
  \item [(a)] $I^M_\mathfrak{p}(N)$ is a $\mathfrak{p}$-secondary submodule of $M$.
  \item [(b)] $Ann_R(I^M_\mathfrak{p}(N))$ is a $\mathfrak{p}$-primary ideal of $R$.
  \item [(c)] $\sqrt{Ann_R(I^M_\mathfrak{p}(N))}=\mathfrak{p}$.
  \item [(d)] $\mathfrak{p}$ is minimal prime ideal of $Ann_R(N)$.
 \end{itemize}
In particular, $I^M_\mathfrak{p}(N)$ is $\mathfrak{p}$-second if and only if $Ann_R(I^M_\mathfrak{p}(N))=\mathfrak{p}$.
\end{thm}

\begin{lem}\label{l3.1}\cite[Lemma 3.1]{MR3588217}
Let $R$ be an integral domain and let $M$ be an Artinian non-zero $R$-module.
\begin{itemize}
  \item [(a)] If $I^M_0(M) \not =0$, then $I^M_0(M)$ is a maximal $(0)$-second submodule of $M$ and it contains every $(0)$-second submodule of $M$.
  \item [(b)] $I^M_0(M) =M$ if and only if $M$ is a $(0)$-second submodule of $M$.
\end{itemize}
\end{lem}

\begin{thm}\label{t3.2}\cite[Theorem 3.2]{MR3588217}
Let $R$ be an integral domain of dimension 1, $M$ be a non-zero Artinian $R$-module and $0 \not = \mathfrak{p} \in V(Ann_R(M))$. Then $I^M_\mathfrak{p}((0:_M\mathfrak{p}))$ is a maximal $\mathfrak{p}$-second submodule of $M$ if and only if $I^M_\mathfrak{p}((0:_M\mathfrak{p})) \not \subseteq I^M_0(M)$.
\end{thm}

\begin{prop}\label{p3.3}\cite[Proposition 3.3]{MR3588217}
Let $Y$ ba a set of prime ideals of $R$ which contains all the maximal ideals, $M$ be an Artinian $R$-module, and $N$ be a non-zero submodule of $M$. Then
$N=\sum_{P \in Y}I^M_P(N)$.
\end{prop}

\begin{cor}\label{c3.4}\cite[Corollary 3.4]{MR3588217}
Let $(R,m)$ be a local ring, $M$ an Artinian $R$-module, and $0 \not = N\leq M$. Then $N$ is cotorsion-free w.r.t. $m$.
\end{cor}

\begin{prop}\label{p3.5}\cite[Proposition 3.5]{MR3588217}
Let $M$ be an Artinian $R$-module, $\mathfrak{p} \in Spec(R)$, and $N$ be a non-zero submodule of $M$.
If $\mathfrak{p}$ is a minimal prime ideal of $Ann_R(N)$ and $I^M_\mathfrak{p}((0:_N\mathfrak{p}))\not=0$, then $I^M_\mathfrak{p}((0:_N\mathfrak{p}))$ is a maximal $\mathfrak{p}$-second submodule of $K\leq M$ such that $I^M_\mathfrak{p}((0:_N\mathfrak{p}))\subseteq K \subseteq N$.
In particular
$I^M_\mathfrak{p}((0:_N\mathfrak{p}))$ is a maximal $\mathfrak{p}$-second submodule of $sec(N)$.
\end{prop}

The following example shows that the condition $I^M_P((0:_NP))\not=0$ in the statement
of Proposition \ref{p3.5} can not be dropped.

\begin{ex}\label{ee3.5}\cite[Example 3.6]{MR3588217}
Consider $M = N = \Bbb Z_{p^{\infty}}$ as $\Bbb Z$-module, where $p$ is a prime number.
Let $q \not= p$ be an another prime number. Then clearly, $q\Bbb Z$ is a minimal prime ideal of
$Ann_{\Bbb Z}(M)$ and $I^M_{(q)}((0 :_N q\Bbb Z)) = (0)$.
\end{ex}

The next theorem gives an important information on the maximal second submodules of $M$.
\begin{thm}\label{t3.6}\cite[Theorem 3.7]{MR3588217}
Let $N$ be a non-zero submodule of an Artinian $R$-module $M$. Then every maximal second submodule of $N$ must be of the form $I^M_P((0:_NP))$ for some $\mathfrak{p} \in V(Ann_R(N)$.
\end{thm}

\begin{defn} \label{d3.17}\cite[Definition 2.11]{AF11}
A second submodule $S$ of an
$R$-module $M$ is said to be a \emph{second sum-irreducible submodule} of $M$ if $S
\not = K+L$, where $K$ and $L$ are second submodules of $M$
properly contained
in $S$.
\end{defn}

\begin{thm} \label{t3.18} \cite[Theorem 2.12]{AF11}
Let $M$ be an Artinian $R$-module. Then every
second submodule of $M$ is a finite sum of second sum-irreducible
submodules of $M$. Moreover, every second sum-irreducible
submodule of $M$ is a sum-irreducible submodule of $M$.
\end{thm}

\begin{defn}\label{d3.10}\cite[Definition 2.8.]{AF12}
Let $R$ be an integral domain. We say that a submodule $N$
an $R$-module $M$ is a \emph{cotorsion-free submodule of $M$} (the dual of torsion-free) if $I^M_0(N)=N$ and is a \emph{cotorsion submodule of $M$} (the dual of torsion) if $I^M_0(N)=0$. Also, $M$ said to be \emph{cotorsion} (resp. \emph{cotorsion-free}) if $M$ is a cotorsion (resp. cotorsion-free) submodule of itself.
\end{defn}

\begin{ex}\label{e3.11}
\begin{itemize}
\item [(a)] The $\Bbb Z$-module $\Bbb Z_{p^\infty}$  (resp. $\Bbb Z$) is
a cotorsion-free (resp. cotorsion) module \cite{AF12}.
\item [(b)] If $R$ is an integral domain, then every non-faithful R-module is cotorsion \cite{AF12}.
\end{itemize}
\end{ex}

An $R$-module $M$ is said to be a
\emph{multiplication module} if for every submodule $N$ of $M$
there exists an ideal $I$ of $R$ such that $N=IM$ \cite{Ba81}.

\begin{prop}\label{p3.12}\cite[Proposition 2.9]{AF12}
Let $R$ be an integral domain and $M$ be a
non-zero $R$-module. Then we have the following.
\begin{itemize}
  \item [(a)] If $M$ is an Artinian faithful multiplication module,
  then $I^M_0(M)$ is a minimal submodule of $M$.
  \item [(b)] If $M$ is Artinian and cotorsion, then
  $Hom_R(S^{- 1}R,M)=0$, where $S=R-0$.
  \item [(c)] If $M$ is an  injective $R$-module and
  $Hom_R(S^{-1}R,M)=0$, where $S=R-0$, then
  $M$ is cotorsion.
  \item [(d)] If $M$ is an Artinian $R$-module, then $M/I^M_0(M)$ is cotorsion.
  \item [(e)] If $M$ is a faithful comultiplication $R$-module, then $M$ is cotorsion-free.
  \end{itemize}
\end{prop}

\begin{cor}\label{t2.4}\cite[Corollary 2.11]{AF12}
Let $R$ be an integral domain and let $M$ be a torsion-free and cotorsion $R$-module.  Then $M$ has no second submodule.
\end{cor}

\begin{ex}\label{t2.5}\cite{AF12}
The $\Bbb Z$-module $\Bbb Z$ is a torsion-free and cotorsion module and it has no second submodule.
\end{ex}
\section{The second socles of submodules}
\noindent
\begin{defn}\label{d3.5}\cite{AF11, CAS103}
 For a submodule $N$ of an $R$-module $M$, the \emph{second radical} (or \emph{second socle}) of $N$ is defined  as the sum of all second submodules of $M$, contained in $N$, and it is denoted by $sec(N)$ (or $soc(N)$). In case $N$ does not contain any second submodule, the second radical of $N$ is defined to be $(0)$. $N \not =0$ is said to be a \emph{second radical} or \emph{socle} submodule of $M$ if $sec(N)=N$.
\end{defn}

\begin{prop}\label{p2.u2}\cite[Proposition 2.1]{AF25}
Let $N$ and $K$ be two submodules of an $R$-module $M$. Then we have the following.
\begin{itemize}
  \item [(a)] If $N\subseteq K$, then $soc(N)\subseteq soc(K)$.
  \item [(b)] $soc(N)\subseteq N$.
  \item [(c)] $soc(soc(N))=soc(N)$.
  \item [(d)] $soc(N)+soc(K)\subseteq soc(N+K)$.
  \item [(e)] $soc(N\cap K)=soc(soc(N) \cap soc(K))$.
  \item [(f)] $soc((0:_MI))=soc((0:_M\sqrt{I}))$ for each ideal $I$ of $R$.
  \item [(g)] $soc((0:_MI^n))=soc((0:_MI))$ for every positive integer $n$ and ideal $I$ of $R$.
  \item [(h)] $soc(N)\subseteq (0:_M\sqrt{Ann_R(N)})$.
  \item [(i)] If $S$ is a $\mathfrak{p}$-second submodule of $M$ such that $S\subseteq N+K$ and $Ann_R(N) \not \subseteq \mathfrak{p}$, then $S \subseteq K$.
 \item [(j)] If $N+K=soc(N)+soc(K)$, then $soc(N+K)=N+K$.
\end{itemize}
\end{prop}

\begin{rem}\label{r2.3}\cite[Remark 2.2]{AF25}
The converse of part (h) of Proposition \ref{p2.u2} is not true in general. For example,
for the submodule $N=0\oplus \Bbb Z_{p^\infty}$ of the $\Bbb Z$-module $M=\Bbb Z_{p^\infty} \oplus \Bbb Z_{p^\infty}$, we have $\sqrt{Ann_{\Bbb Z}(N)}=0$ (here $\Bbb Z$ denotes the ring of integers). Thus $N=soc(N) \not = (0:_M\sqrt{Ann_{\Bbb Z}(N)})$.
\end{rem}

An $R$-module $M$ is said to be \emph{atomic} if every nonzero submodule of $M$ contains a
minimal submodule \cite{HP10}.
\begin{thm}\label{c2.4}\cite[Theorem 2.3]{AF25}
Let $M$ be an atomic $R$-module. Then we have the following.
\begin{itemize}
\item [(a)] $soc(M)=0$ if and only if $M=0$.
\item [(b)] If $N$ and $K$ are two submodules of $M$,  then $soc(N) \cap soc(K)=0$ if and only if $N \cap K=0$.
\item [(c)] If $m$ a maximal ideal of $R$ and $Q$ is an $m$-secondary submodule of $M$, then $soc(Q)$ is an $m$-second submodule of $M$.
\item [(d)] If $R$ is an integral domain with $dim R=1$ and $M$ is a secondary $R$-module, then $soc(M)$ is a second submodule of $M$.
\end{itemize}
\end{thm}

\begin{prop}\label{p2.6}\cite[Proposition 2.4]{AF25}
Let $M$ be an $R$-module. Then the following hold.
\begin{itemize}
  \item [(a)] If $R$ is an Artinian ring, then $Soc_R(M)=soc(M)$.
  \item [(b)] If $N$ is a submodule of $M$, then $soc(N)$ is the sum of the maximal second submodules of $N$.
  \item [(c)] If $Q$ is a $\mathfrak{p}$-secondary submodule of $M$, then $soc(Q)=soc(Q \cap (0:_M\mathfrak{p}))$.
  \item [(d)] If $(R,m)$ is a local ring, then every second submodule of $M$ is minimal if and only if every socle submodule of $M$ is minimal.
  \item [(e)] If $(R,m)$ is a local ring and every second submodule of $M$ is minimal, then $soc(M)$ is simple or zero.
\end{itemize}
\end{prop}

\begin{prop}\label{p2.7}\cite[Proposition 2.5]{AF25}
If $V$ is a vector space, then $soc(N_1
+N_2)=soc(N_1)+soc(N_2)$ for every pair of subspaces $N_1$, $N_2$ of $V$.
\end{prop}

\begin{thm}\label{t2.9}\cite[Theorem 2.6]{AF25}
Let $M$ be an $R$-module. Then the following hold.
\begin{itemize}
  \item [(a)] Let $I$ be an ideal of $R$ and $N$ be a submodule of $M$. If $S$ is a $\mathfrak{p}$-second submodule of $M$ such that $S\subseteq (0:_MI)+N$, then
$S\subseteq (0:_MI)$ or $S\subseteq N$.
  \item [(b)] Let $N$ and $K$ be two submodules of $M$ such that whenever
$S\subseteq N+K$, we have $S \subseteq N$ or $S\subseteq K$
for every second submodule $S$ of $M$. Then
$$
soc(N+K)=soc(N)+soc(K).
$$
\item [(c)] Let $M$ be a comultiplication $R$-module. If $S$ is a second submodule of $M$ such that $S\subseteq N+K$ for any pair of submodules $N$ and $K$ of $M$, then either $S\subseteq N$ or $S\subseteq K$. Consequently,
      $$
      soc(N+K)=soc(N)+soc(K)
      $$
  for every pair of submodules $N$ and $K$ of $M$.
\item [(d)] $soc((0:_MI)+N)=soc((0:_MI))+soc(N)$
for every ideal $I$ of $R$ and every submodule $N$ of $M$.
\end{itemize}
\end{thm}

\begin{thm}\label{t2.13}\cite[Theorem 2.7]{AF25}
Let $M$ be an $R$-module. If $N$ and $K$ are two submodules of $M$ such that $Ann_R(N)$ and $Ann_R(K)$ are comaximal, then $soc(N+K)=soc(N)+soc(K)$.
\end{thm}

\begin{cor}\label{c2.14}\cite[Corollary 2.8]{AF25}
Let $K_1, ..., K_n$ be submodules of an $R$-module $M$ such that $Ann_R(K_i)$ are pairwise comaximal. Then $soc(K_1+...+K_n)=soc(K_1)+...+soc(K_n)$.
\end{cor}

\begin{defn}\label{c2.1uiy4}\cite{AF21}
Let $M$ be an $R$-module. The set of all second submodules of $M$ is
called the \emph{second spectrum} of $M$ and denoted by $Spec^s(M)$. The map $\phi: Spec^s(M)\rightarrow Spec(R/Ann_R(M))$ defined by $\phi (S)=Ann_R(S)/Ann_R(M)$ for every $S \in Spec^s(M)$, is called the \emph{natural map} of $Spec^s(M)$.
\end{defn}

\begin{prop}\label{p2.15}\cite[Proposition 2.9]{AF25}
Let $M$ be an $R$-module and let $N$ be a submodule of $M$ such that the natural map $\phi$ of $Spec^s(N)$ is surjective. Then $Ann_R(soc(N))=\sqrt{Ann_R(N)}$.
\end{prop}

\begin{lem}\label{l2.16}\cite[Lemma 2.10]{AF25}
Let $M$ be an $R$-module. Then we have the following.
\begin{itemize}
   \item [(a)] If $M$ is a finitely generated comultiplication module
and $\mathfrak{p}$ is a prime ideal of $R$ containing $Ann_R(M)$, then $(0:_M\mathfrak{p})$ is a second
submodule of $M$.
 \item [(b)] If $M$ is a finitely generated comultiplication module, then the natural map $\phi$ of $Spec^s(M)$ is surjective.
  \item [(c)] If the natural map $\phi$ of $Spec^s(M)$ is surjective and $I$ is an ideal of $R$ containing $Ann_R(M)$, then $Ann_R((0:_M\sqrt{I}))=\sqrt{I}$.
\end{itemize}
\end{lem}

\begin{thm}\label{t2.17}\cite[Theorem 2.11]{AF25}
Let $M$ be a faithful $R$-module such that the natural map $\phi$ of $Spec^s(M)$ is surjective.
Consider the following equalities:
\begin{itemize}
  \item [(a)] $soc((0:_MI))=(0:_M\sqrt{I})$ for each ideal $I$ of $R$.
  \item [(b)] $soc(N)=(0:_M\sqrt{Ann_R(N)})$ for each submodule $N$ of $M$.
  \item [(c)] $Ann_R(soc(N))=\sqrt{Ann_R(N)}$ for each submodule $N$ of $M$.
  \item [(d)] $Ann_R(soc((0:_MI)))=\sqrt{I}$ for each ideal $I$ of $R$.
\end{itemize}
Then $(b)\Rightarrow (c)\Rightarrow (d)$ and $(b)\Rightarrow (a)\Rightarrow (d)$. Furthermore, if $M$ is a comultiplication module, then (a), (b), (c) and (d) are all equivalent.
\end{thm}

\begin{thm}\label{t2.11}\cite[Theorem 2.12]{AF25}
Let $N$ and $K$ be two submodules of a finitely generated comultiplication $R$-module $M$. Then the following hold.
\begin{itemize}
  \item [(a)] $soc(N)=(0:_M\sqrt{Ann_R(N)})$.
  \item [(b)] $Ann_R(soc(N))=\sqrt{Ann_R(N)}$.
  \item [(c)] If $Ann_R(K)=\sqrt{Ann_R(K)}$ and $Ann_R(N)=\sqrt{Ann_R(N)}$, then $Ann(soc(N+K))=Ann_R(N+K)$.
  \item [(d)] If $N$, $K$ are secondary submodule of $M$ with $soc(N)=soc(K)$, then $N+K$ is a secondary submodule of $M$.
\end{itemize}
\end{thm}

\begin{cor}\label{c2.12}\cite[Corollary 2.13]{AF25}
If $Q$ is a secondary submodule of a finitely generated comultiplication $R$-module $M$,
then $soc(Q)$ is a second submodule of $M$.
\end{cor}

Let $N$ and $K$ be two submodules of an $R$-module $M$. In \cite{Lu90}, it is shown that in general, $rad(N \cap K) \not = rad(N) \cap rad(K)$. We haven't found any example of an $R$-module $M$ such that for some submodules $N$ and $K$ of $M$, $soc(N+K)\not=soc(N)+soc(K)$. This motivates the following question.

\begin{co}\label{lq2.20}\cite[Question 2.14]{AF25}
Let $N$ and $K$ be two submodules of an $R$-module $M$. Is $soc(N+K)=soc(N)+soc(K)$?
\end{co}

\begin{lem}\label{l2.20}\cite[Lemma 3.1]{AF25}
If $R$ is an integral domain and $M/I^M_0(M)$ is a
finitely cogenerated $R$-module such that $I^M_0(M) \not= 0$, then $I^M_0(M)$ is a (0)-second submodule of $M$.
\end{lem}

\begin{thm}\label{t2.21}\cite[Theorem 3.2]{AF25}
Let $R$ be a Noetherian integral domain and $M$ be a finitely cogenerated $R$-module such that $I^M_0(M) \not =0$. If $I^M_0(M)$ is a finitely generated $R$-module with only finitely many second submodules, then for any secondary submodule $Q$ of $M$, $soc(Q)$ is second.
\end{thm}

\begin{defn}\label{d2.18}\cite[Definition 3.3]{AF25}
The \emph{second submodule dimension} of an $R$-module $M$, denoted by $S.dim M$, is defined to be the supremum of the length of chains of second submodules of $M$ if $Spec^s(M) \not =\emptyset$ and $-1$ otherwise.
\end{defn}

\begin{thm}\label{t2.188}\cite[Theorem 3.4]{AF25}
Let $K$ be a field and $M$ a $K$-Vector space with $dim_KM=n$. Then $S.dim M=n-1$.
\end{thm}

\begin{thm}\label{t2.19}\cite[Theorem 3.5]{AF25}
If $R$ is a one dimensional domain and $M$ is a finitely cogenerated cotorsion module with
$S.dim M=1$, then the following are equivalent.
\begin{itemize}
  \item [(a)] $M$ is a second module.
  \item [(b)] $S_1+S_2=M$ for any distinct second submodules $S_1$ and $S_2$.
  \item [(c)] Every proper submodule contains exactly one second submodule.
  \item [(d)] Every proper second submodule is minimal.
\end{itemize}
\end{thm}

\begin{cor}\label{c2.119}\cite[Corollary 3.6]{AF25}
If $R$ is a one dimensional domain and $M$ is a second finitely cogenerated cotorsion module with
$S.dim M=1$, then $soc(N)$ is second for every non-zero submodule $N$ of $M$.
\end{cor}

\begin{thm}\label{t2.22}\cite[Theorem 3.7]{AF25}
Let $R$ be an integral domain and $M$ an $R$-module with $S.dim M=1$. Then the following hold.
\begin{itemize}
  \item [(a)] If $M$ is Artinian cotorsion-free, then $soc(N)=Soc_R(N)$ for any proper submodule $N$ of $M$.
  \item [(b)] If $R$ is  Noetherian one-dimensional and $M$ is a second finitely cogenerated cotorsion module, then every non-zero submodule $N$ of $M$ with $\sqrt{Ann_R(N)}\not =0$ is secondary.
\end{itemize}
\end{thm}

\begin{thm}\label{t3.6}\cite[Theorem 2.5]{AF11}
 Let $M$ be an $R$-module. If $M$ satisfies
the descending chain condition on socle submodules, then every
non-zero submodule of $M$ has only a finite number of maximal
second submodules.
\end{thm}

\begin{cor}\label{c3.7}\cite[Corollary 2.6]{AF11}
 Every Artinian $R$-module
contains only a finite number of maximal second submodules.
\end{cor}

\begin{thm} \label{t3.7}\cite[Theorem 2.1]{AF12}
 Let $R$ be a Noetherian domain, $M$ an Artinian $R$-module, and $N$ a proper submodule of $M$. If $Ann_R(N)$ is a radical ideal contained in only finitely many prime ideals, each of which is maximal, then $N$ is a socle submodule of $M$.
\end{thm}

\begin{lem}\label{l2.9}\cite[Lemma 2.2]{AF12}
Let $M$ be a finitely generated second $R$-module. Then $M$ is an Artinian
$R$-module.
\end{lem}

\begin{thm} \label{t3.8}\cite[Theorem 2.3]{AF12}
Let $M$ be a Noetherian $R$-module.
Then $M$ satisfies the descending chain condition on socle submodules.
\end{thm}

\begin{cor}\label{c3.9}\cite[Corollary 2.4]{AF12}
Every Noetherian $R$-module has only a finite number of maximal second submodules.
\end{cor}

\begin{thm}\label{t3.2}\cite[Theorem 2.5]{AF12}
Let $M$ be a faithful finitely generated comultiplication $R$-module
satisfying the descending chain condition on second submodules. Then $R$ satisfies
the ascending chain condition on prime ideals..
\end{thm}

\begin{prop}\label{p2.7}\cite[Proposition 2.7]{AF12}
Let $M$ be an $R$-module. Then we have the following.
\begin{itemize}
  \item [(a)] If $M$ is Artinian, then there exist $n \in \Bbb N$ and prime ideals $\mathfrak{p}_i$, $1\leq i \leq n$, such that $Ann_R((0:_M\mathfrak{p}_i))=\mathfrak{p}_i$ and $soc(M)=\sum^n_{i=1}I^M_\mathfrak{p}((0:_M\mathfrak{p}_i))$.
  \item [(b)] If $R$ is a one-dimension integral domain, then $soc(M)\subseteq I^M_0(M)+Soc_R(M)$. Moreover, if $M/I^M_0(M)$ is finitely cogenerated the equality holds.
  \item [(c)] If $\mathfrak{p}$ is a prime ideal of $R$ and $M$ has a finitely generated $\mathfrak{p}$-second submodule, then $I^M_\mathfrak{p}(0:_M\mathfrak{p})=(0:_M\mathfrak{p})$.
\end{itemize}
 \end{prop}

\begin{thm}\label{t5.11}\cite[Theorem 2.21]{MR3755273}
Let $M$ be a finitely generated comultiplication $R$-module and $N$ be a submodule of $M$. Then $sec(M) \subseteq N$ if and only if $Ann_R(N) \subseteq \sqrt{Ann_R(M/N)}$.
\end{thm}

\begin{prop}\cite[Proposition 3.1]{MR3307390}
Let $M$ be a faithful Noetherian
comultiplication $R$-module, then $sec(M) = soc(M)=(0:_M Rad (R))$, where
$Rad(R)$ is the Jacobson radical of $R$.
\end{prop}

\begin{cor}\cite[Corollary 3.13]{MR3307390}
Let $R$ be a Dedektnd domatn and $M$ be an amply supplemented
second radical $R$-module with finite hollow dimension. Then $M$ is a
finite sum of minimal second submodules.
\end{cor}

\begin{cor}\label{c3.7}\cite[Corollary 3.8]{MR3588217}
Let $M$ be an Artinian $R$-module and $0 \not = N\leq M$. Then $soc(N)=\sum_{\mathfrak{p} \in Y}I^M_\mathfrak{p}((0:_N\mathfrak{p}))$, where $Y$ is a finite subset of $V(Ann_R(N))$.
\end{cor}

\begin{cor}\label{c3.8}\cite[Corollary 3.9]{MR3588217}
Let $N$ be a non-zero submodule of an Artinian $R$-module $M$. If $N$ is a $\mathfrak{p}$-secondary submodule of an $R$-module $M$ for some $\mathfrak{p} \in Spec(R)$, then we have the following.
\begin{itemize}
  \item [(a)] $I^M_\mathfrak{p}((0:_N\mathfrak{p}))$ is a maximal $\mathfrak{p}$-second submodule of $soc(N)$.
  \item [(b)] If $\mathfrak{p}$ is a maximal ideal of $R$, then $soc(N)=I^M_\mathfrak{p}((0:_N\mathfrak{p}))$ so that $soc(N)$ is a $\mathfrak{p}$-second submodule of $M$.
\end{itemize}
\end{cor}

\begin{cor}\label{c3.9}\cite[Corollary 3.10]{MR3588217}
Let $I$ be an ideal of a ring $R$ and $M$ be an Artinian $R$-module such that $(0:_MI) \not =0$. Then
$$
 soc((0:_MI))=\sum_{\mathfrak{p} \in V(Ann_R((0:_MI))}I^M_\mathfrak{p}((0:_M\mathfrak{p})).
$$
\end{cor}

\begin{ex}\label{e3.10}\cite[Example 3.11]{MR3588217}
For any prime integer $p$, let $M=(\Bbb Z/p\Bbb Z)\times \Bbb Z_{p^\infty}$. Then $M$ is an Artinian faithful $\Bbb Z$-module and $V(Ann_{\Bbb Z}(M))=V(0)=Spec(\Bbb Z)$. Hence $soc(M)=\sum_{(q) \in V(0)}I^M_{(q)}((0:_Mq))$ by Corollary \ref{c3.9}. Since $I^M_{(q)}((0:_Mq))=I^M_{(q)}(0)=0$ for each $q \not = p$,
$$
soc(M)=I^M_{(0)}(M)+I^M_{(p)}((0:_MP))=
$$
$$
((\Bbb Z/p\Bbb Z)\times \Bbb Z_{p^\infty})+((\Bbb Z/p\Bbb Z)\times <1/p+\Bbb Z>)=M.
$$
\end{ex}

\begin{thm}\label{c103}\cite[Corollary 3.5]{CAS103}
Let $R$ be a commutative one-dimensional Noetherian domain and $M$ be an $R$-module. Then $sec(M)=Soc(M)+div(M)=K\oplus div(M)$ for some semisimple submodule $K$ of $M$, where $div(M)$ denotes the sum of all divisible submodules of $M$.
\end{thm}

\noindent
 A proper submodule $L$ of an $R$-module $M$ is \emph{radical} if $L$ is an intersection of prime submodules of $M$. Moreover, a submodule $L$ of $M$ is \emph{isolated} if, for each proper submodule $N$ of $L$, there exists a prime submodule $K$ of $M$ such that $N\subseteq K$ but $L \not \subseteq K$ \cite{MS06}.
\begin{defn}\label{d4.1} \cite[Definition 3.1]{AF11}
We say that a submodule $N$ of an
$R$-module $M$ is \emph{coisolated} if $soc(L) \not = soc(N)$
for every submodule $L$ of $M$ that properly contains $N$. This can be regarded as a dual notion of the isolated submodule.
\end{defn}

\begin{exs}\label{e4.2}\cite[Examples 3.2]{AF11}
\begin{itemize}
   \item [(a)] Every submodule of the $\Bbb Z$-module $\Bbb
Z_6$ is coisolated.
   \item [(b)]  Every submodule of the $\Bbb Z$-module $\Bbb
Z_{p^\infty}$ is not coisolated.
   \item [(c)] $\Bbb Z_{p^\infty}$ as a $\Bbb Z$-module is a coisolated submodule of  $\Bbb Z_{p^\infty}$ but it is not an isolated submodule of $\Bbb
       Z_{p^\infty}$.
   \item [(d)] $(0)$ is an isolated submodule of $\Bbb Z$ but it is not a coisolated submodule of $\Bbb Z$.
\end{itemize}
\end{exs}

\begin{defn}\label{d4.7}\cite[Definition 3.3]{AF11}
Let $N$  and $K$ be two submodules of an $R$-module $M$
with $N \subset K$ and let $K/N$ be a second submodule of $M/N$. Then we
say $K/N$ can be \emph{lowered} to $M$ if there exists a
second submodule $L$ of $M$
such that $K=L+N$. This notion can be
regarded as a dual notion of lifting \cite{MS06}.
\end{defn}

\begin{thm} \label{t4.8}\cite[Theorem 3.4]{AF11}
(Dual of Lying over.) Let $\mathfrak{p}$ be a prime ideal
of $R$ and $M$ be an $R$-module which satisfies the property $AB5^*$.
Let $N$ and $K$ be two submodules of $M$ with $N\subset K$.
Then a $\mathfrak{p}$-second submodule $K/N$ of $M/N$ can be lowered to $M$ if and
only if $K\subseteq N+(0:_M\mathfrak{p})$.
\end{thm}

\begin{prop}\label{t4.9} \cite[Proposition 3.5]{AF11}
A submodule $N$ of an $R$-module $M$ is
coisolated if and only if for each submodule $H$ of $M$ with $N
\subset H$ there exists a second submodule $K/N$ of $M/N$ such
that $K\subseteq H$ and $K/N$ can be lowered to $M$.
\end{prop}

A submodule $N$ of an $R$-module $M$ is said to be  \emph{copure} if $(N:_MI)=N+(0:_MI)$ for every ideal $I$ of $R$ \cite{AF09}.
\begin{thm} \label{t4.12}\cite[Theorem 3.6]{AF11}
Let $M$ be an $R$-module satisfying the
property $AB5^*$. Then the following statements are equivalent.
\begin{itemize}
  \item [(a)] Every submodule of $M$ is a coisolated submodule of $M$.
  \item [(b)] Every non-zero submodule of $M$ is a socle submodule of $M$.
  \item [(c)] Every submodule of $M$ is copure.
  \item [(d)] $N+(0:_M\mathfrak{p})=(N:_M\mathfrak{p})$ for every submodule $N$ of $M$
  and every maximal ideal $\mathfrak{p}$ of $R$.
\end{itemize}
\end{thm}

\begin{thm} \label{t4.14}\cite[Theorem 3.7]{AF11}
The following statements are equivalent for a
submodule $N$ of an $R$-module $M$ such that $M/N$ is finitely
cogenerated.
\begin{itemize}
  \item [(a)] $N$ is coisolated.
  \item [(b)] $N$ is a direct summand of $K$ for every submodule $K/N$
  of $M/N$.
  \item [(c)] For every minimal submodule $K/N$ of $M/N$, there exists a
  minimal submodule $L$ of $M$ such that $K=L+N$.
  \item [(d)] $(N:_M\mathfrak{p})=N+(0:_M\mathfrak{p})$ for every maximal ideal $\mathfrak{p}$ of
  $R$.
  \item [(e)] $N$ is a direct summand of
  $(N:_M\mathfrak{p})$ for every maximal ideal $\mathfrak{p}$ of $R$.
\end{itemize}
\end{thm}

\begin{thm} \label{t4.15}\cite[Theorem 3.8]{AF11}
Let $R$ be a PID and let $N$ be a socle submodule of an $R$-module $M$ such that $M/N$ is a finitely cogenerated $R$-module. Then $N$ is coisolated if and only if $N$ is a copure submodule of $M$.
\end{thm}

Recall that a submodule $N$ of an $R$-module $M$ is said to be \emph{pure} if $IN=N \cap IM$ for every ideal $I$ of $R$ \cite{AF74}.
\begin{cor}\label{c4.16}\cite[Corollary 3.9]{AF11}
Let $R$ be a PID and $N$ a finitely generated submodule of an $R$-module $M$ such that $M/N$ is a finitely cogenerated $R$-module. If  $N$ is a socle and a radical submodule of $M$, then the following statements are equivalent.
\begin{itemize}
  \item [(a)] $N$ is an isolated submodule of $M$.
  \item [(b)] $N$ is a coisolated submodule of $M$.
  \item [(c)] $N$ is a pure submodule of $M$.
  \item [(d)] $N$ is a copure submodule of $M$.
\end{itemize}
\end{cor}

\begin{ex}\cite[Example 3.10]{AF11}
Every submodule of the $\Bbb Z$-module $\Bbb Z_k$ is both isolated and coisolated if $k$ is square-free.
\end{ex}
\section{Finitely generated coreduced comultiplication modules}
Let $M$ be an $R$-module.
The set of all maximal second submodules of $M$ will be denoted by $Max^s(M)$. The sum of all maximal second submodules of $M$ contained in a submodule $K$ of $M$ is denote by $\mathfrak{S}_K$.  In case $M$ does not contain any maximal second submodule which is contained in $K$, then $\mathfrak{S}_K$ is defined to be $(0)$. If $N$ is a submodule of $M$, define
 $V^s(N) = \{S \in Max^s(M) : S \subseteq N\}$.
$M$ is said to be \emph{coreduced module} if $(L:_Mr)=M$ implies that $L+(0:_Mr)=M$, where $r \in R$ and $L$ is a completely irreducible submodule of $M$ \cite{MR3755273}.

\begin{prop}\label{p5.12}\cite[Proposition 2.22]{MR3755273}
Let $M$ be an $R$-module. Then $M$ is a coreduced $R$-module if $sec(M)=M$. The converse holds when $M$ is a finitely generated comultiplication $R$-module.
\end{prop}

In \cite{Re03}, Redmond introduced the definition of the zero-divisor graph with respect to an ideal. Let $I$ be an ideal of $R$. The zero-divisor graph of $R$ with respect to $I$, denoted by $\Gamma_I (R)$, is the graph whose vertices are the set $$
\{x \in R \setminus I \mid xy \in I \ for\ some \ y \in R \setminus I \}
$$
with distinct vertices $x$ and $y$ are adjacent if and only if $xy \in I$.

\begin{thm}\label{t5.13}\cite[Theorem 2.23]{MR3755273}
Let $M$ be a finitely generated comultiplication $R$-module and $sec(M) \subseteq N \not = M$. If $\Gamma_{Ann_R(M)}(R)$ is complemented, then there exists $a \in Ann_R(N)$ such that $a^tM=0$, $a^{t-i}M \not =0$ and $a^{t-1} \perp a^i$, $t=2,3$ and $1 \leq i \leq t-2$.
\end{thm}

\begin{lem}\label{l5.14}\cite[Lemma 2.24]{MR3755273}
Let $M$ be a coreduced comultiplication $R$-module and $I$ be an ideal of $R$. If $I \subseteq \mathfrak{p}$, where $\mathfrak{p}$ is a minimal prime ideal of $Ann_R(M)$, Then $I \subseteq W_R(M)$.
\end{lem}

\begin{thm}\label{t5.18}\cite[Theorem 2.25]{MR3755273}
Let $M$ be a finitely generated comultiplication $R$-module. Then we have the following.
\begin{itemize}
  \item [(a)] If $R$ is a ring with $|\bar{R}|> 4$ and $\Gamma_{Ann_R(M)}(R)$ is a complete graph, then either $(0:_MZ_R(\bar{R}))=0$ or
$(0:_MZ_R(\bar{R}))=sec(M)$.
  \item [(b)] If $sec(M) \not =M$ and there are $\alpha , \beta \in V(\Gamma_{Ann_R(M)}(R))$ such that $R\alpha + R\beta \not \subseteq W_R(M)$, then  $diam(\Gamma_{Ann_R(M)}(R))=3$.
\end{itemize}
\end{thm}

\begin{prop}\label{l0.1}\cite[Proposition 2.1]{FF2024}
Let $M$ be a finitely generated comultiplication $R$-module. Then we have the following.
\begin{itemize}
\item [(a)] If $S$ is a maximal second submodule of $N$, then $Ann_R(S)$ is a prime ideal minimal over
$Ann_R(N)$.
\item [(b)] If $S$ is a submodule of $M$ such that $Ann_R(S)$ is a prime ideal minimal over
$Ann_R(M)$, then  $S$ is a maximal second submodule of $M$.
\item [(c)] If $M$ is a coreduced $R$-module, $S$ a maximal second submodule of $M$, and $a \in Ann_R(S)$, then
$$
Ann_{R/Ann_R(M)}(a+Ann_R(M)) \not \subseteq Ann_R(S)/Ann_R(M).
$$
\end{itemize}
\end{prop}

\begin{thm}\label{tt1.11}\cite[Theorem 2.2]{FF2024}
Let $M$ be a faithful finitely generated coreduced comultiplication $R$-module. Then for each $a \in R$,  we have $V^s((0:_MAnn_R(a))=Max^s(M) \setminus V^s((0:_Ma))$.
\end{thm}

The intersection of all minimal prime ideals of $R$ containing an ideal $I$ of $R$ is denote by $\mathfrak{p}_I$.
\begin{thm}\label{t119.3}\cite[Theorem 2.3]{FF2024}
Let $M$ be a faithful finitely generated comultiplication $R$-module. Then $(0:_M\mathfrak{p}_I)=\mathfrak{S}_{(0:_MI)}$ for each ideal $I$ of $R$.
\end{thm}

\begin{thm}\label{t1.5}\cite[Theorem 2.5]{FF2024}
Let $M$ be a finitely generated coreduced comultiplication $R$-module. Then we have the following.
\begin{itemize}
\item [(a)] If $M$ is a faithful $R$-module, then $V^s((0:_Ma))=V^s(Ann_R(a)M)$ for each  $a \in R$.
\item [(b)] $Ann_R(IJM)M=Ann_R(IM)M + Ann_R(JM)M$ for each ideals $I, J$ of $R$.
\end{itemize}
\end{thm}

\begin{prop}\label{p0.1}\cite[Proposition 2.8]{FF2024}
Let $S$ be a second submodule of an $R$-module $M$. Then $S \subseteq I^M_{Ann_R(S)}(M)$.
\end{prop}

For each prime ideal $\mathfrak{p} $ of $R$, set $nil \mathfrak{p}  = \cap \mathfrak{p}'$, where $\mathfrak{p}'$ ranges over all prime ideals of $R$ contained in $\mathfrak{p}$.

Let $\mathfrak{p} $ be a prime ideal of $R$. Then the set
$O_{\mathfrak{p}} = \{a \in R : Ann_R(a)\not \subseteq \mathfrak{p} \}$
is an ideal of $R$ contained in $\mathfrak{p}$ and $\sqrt{O_\mathfrak{p}} = nil \mathfrak{p}$, so if $\mathfrak{p} $ is a minimal
prime ideal of $R$, $\mathfrak{p} = \sqrt{O_\mathfrak{p} }$, in particular, $\mathfrak{p}= O\mathfrak{p}$, when $R$ is a reduced ring \cite{MR2839935}.
\begin{note}\cite[Notation and Remark 2.9]{FF2024}
Let $S$ be a second submodule of an $R$-module $M$.
We define $conil(S)=\sum S' $, where $S' $ ranges over all second submodules of $M$ such that $S\subseteq S' $.
\end{note}

We set $\overline{I}=I+Ann_R(M)$ for each ideal $I$ of $R$ and $\overline{a}=a+Ann_R(M)$ for each $a \in R$.
\begin{thm}\label{t0.2}\cite[Theorem 2.10]{FF2024}
Let $S$ be a second submodule of a finitely generated comultiplication $R$-module $M$. Then we have the following.
\begin{itemize}
\item [(a)]
$$
conil(S) = (0:_Mnil(\overline{Ann_R(S)}) =(0:_M\sqrt{O_{\overline{Ann_R(S)}}}).
$$
\item [(b)]
$$
sec(I^M_{Ann_R(S)}(M))=conil(S)=(0:_M\sqrt{O_{\overline{Ann_R(S)}}}).
$$
In particular, if $M$ is coreduced, then $sec(I^M_{Ann_R(S)}(M))=I^M_{Ann_R(S)}(M)$.
\end{itemize}
\end{thm}

\begin{cor}\label{c00.3}\cite[Corollary 2.11]{FF2024}
Let $S$ be a submodule of a finitely generated comultiplication $R$-module $M$. Then $S$ is a maximal second submodule of $M$ if and only if
$S=sec(I^M_{Ann_R(S)}(M))$. In particular, if $M$ is a coreduced $R$-module, then $S$ is a maximal second submodule of $M$
if and only if $S=I^M_{Ann_R(S)}(M)$.
\end{cor}

\begin{note}\cite[Notation and Remark 2.12]{FF2024}
For a submodule $N$ of an $R$-module $M$, we define
$$
V^*(N) = \{S \in Max^s(M) : S \subseteq N \ and \ (N:_RM) \not \subseteq Ann_R(S)\}.
$$
Clearly, $V^*(N)\subseteq V^s(N)$.
\end{note}

\begin{thm}\label{t0.3}\cite[Theorem 2.13]{FF2024}
Let $L$ be a completely irreducible submodule of a comultiplication $R$-module $M$. Then $V^s(L) \subseteq V^*(L)$  if one of the following conditions hold.
\begin{itemize}
\item [(a)] $M$ is a faithful finitely generated coreduced $R$-module and $Ann_R(L)$ is a finitely generated ideal of $R$.
\item [(b)] $M$ is a finitely generated coreduced $R$-module and $M/S$ is a finitely cogenerated $R$-module for each maximal second submodule $S$ of $M$ with $S \subseteq L$.
\item [(c)] $L$ is a pure submodule of $M$.
\end{itemize}
 \end{thm}

\begin{prop}\label{lll0.3}\cite[Proposition 2.14]{FF2024}
Let $M$ be an $R$-module such that $V^s(L) \subseteq V^*(L)$ for each completely irreducible submodule $L$ of $M$. If $K$ be a submodule of $M$ such that $M/K$ is a finitely cogenerated $R$-module and $S$ is a maximal second submodule of $M$ with $S \subseteq K$, then $(K:_RM) \not \subseteq Ann_R(S)$.
 \end{prop}

\begin{cor}\label{cl0.3}\cite[Corollary 2.15]{FF2024}
Let $M$ be a Noetherian coreduced comultiplication $R$-module. Then $V^*(K)= V^s(K)$ for each submodule $K$ of $M$.
 \end{cor}

\begin{thm}\label{tt0.1}\cite[Theorem 2.17]{FF2024}
Let $M$ be a finitely generated faithful coreduced comultiplication $R$-module. Then $\mathfrak{S}_{IM}= sec(IM)=IM$ for each ideal $I$ of $R$.
\end{thm}

\begin{cor}\label{ccc1.13}\cite[Corollary 2.18]{FF2024}
Let $M$ be a faithful finitely generated coreduced comultiplication and multiplication $R$-module. Then for each submodule $N$ of $M$, we have  $\mathfrak{S}_{N}= sec(N)=N$.
\end{cor}

\begin{ex}\label{e3.9}\cite[Example 2.19]{FF2024}
Consider the $\Bbb Z_n$-module $M=\Bbb Z_n$, where $n$ is square free.
We know that $M$ is a faithful finitely generated coreduced comultiplication and multiplication $\Bbb Z_n$-module. Thus
for each submodule $N$ of $M$, we have  $\mathfrak{S}_{N}= sec(N)=N$.
\end{ex}

\begin{cor}\label{c1.13}\cite[Corollary 2.20]{FF2024}
Let $M$ be a faithful finitely generated coreduced comultiplication $R$-module. Then for each $a \in R$, $Ann_R(aM)M=\mathfrak{S}_{(0:_Ma)}$.
\end{cor}

The following example shows that the condition "$M$ is a finitely generated $R$-module" in Corollary \ref{c1.13} is necessary.
\begin{ex}\label{e1.13}\cite[Example 2.21]{FF2024}
For each prime number $p$ the $\Bbb Z$-module $\Bbb Z_{p^\infty}$ is a faithful coreduced comultiplication $\Bbb Z$-module. But
the $\Bbb Z$-module $\Bbb Z_{p^\infty}$ is not finitely generated.
For each positive integer $n$,
$$
Ann_{\Bbb Z}(p^n\Bbb Z_{p^\infty})\Bbb Z_{p^\infty}=0\not =\langle1/p+ \Bbb Z\rangle=\mathfrak{S}_{(0:_{\Bbb Z_{p^\infty}}p^n)}.
$$
\end{ex}

\begin{thm}\label{t1.45}\cite[Theorem 2.22]{FF2024}
Let $M$ be a finitely generated coreduced comultiplication $R$-module. Then for each submodule $K$ of $M$ we have the following.
\begin{itemize}
\item [(a)] $V^s((0:_M(K:_RM))=Max^s(M) \setminus V^*(K)$.
\item [(b)] $V^*(K)=V^*((K:_RM)M)$.
\end{itemize}
\end{thm}

\begin{cor}\label{c1.143}\cite[Corollary 2.23]{FF2024}
Let $M$ be a Noetherian coreduced comultiplication $R$-module. Then for each submodule $K$ of $M$, $(K:_RM)M=\mathfrak{S}_{K}$.
\end{cor}

\begin{thm}\label{t99.6}\cite[Theorem 2.24]{FF2024}
Let $N$ be a submodule of a Noetherian coreduced comultiplication $R$-module $M$. Then for submodules $K$ and $H$ of $M$ the following are equivalent:
\begin{itemize}
\item [(a)]  $\mathfrak{S}_K=\mathfrak{S}_H$ and $N \subseteq  K$ imply that $N  \subseteq H$;
\item [(b)] $V^s(K)=V^s(H)$ and $N \subseteq K$ imply that $N  \subseteq H$;
\item [(c)]  $(K:_RM)=(H:_RM)$ and $N  \subseteq  K$ imply that $N  \subseteq H$;
\item [(d)] For submodule $K$ of $M$, $N \subseteq K$ implies that $N \subseteq (K:_RM)M$.
\end{itemize}
\end{thm}

\section{$S$-second submodules of a module}
\begin{thm}\label{t2.333}\cite[Theorem 2.2]{MR4401391}
Let $S$ be a m.c.s. of $R$.  For a submodule $N$ of an $R$-module $M$  with $Ann_R(N) \cap S= \emptyset$ the following statements are equivalent:
\begin{itemize}
 \item [(a)]  There exists an $s \in S$ such that $srN=sN$ or $srN=0$ for each $r \in R$;
 \item [(b)]  There exists an $s \in S$ and whenever $rN\subseteq K$, where $r \in R$ and $K$ is a submodule of $M$, implies either that $rsN=0$ or $sN\subseteq K$;
 \item [(c)] There exists an $s \in S$ and whenever $rN\subseteq L$, where $r \in R$ and $L$ is a
completely irreducible submodule of $M$, implies either that $rsN=0$
 or $sN\subseteq L$.
 \item [(d)] There exists an $s \in S$, and $JN\subseteq K$ implies $sJ \subseteq Ann_R(N)$ or $sN \subseteq K$ for each ideal $J$ of $R$ and submodule $K$ of $M$.
\end{itemize}
\end{thm}

\begin{defn}\label{d2.1}\cite[Definition 2.3]{MR4401391}
Let $S$ be a m.c.s. of $R$ and $N$ be a submodule of an $R$-module $M$ such that $Ann_R(N) \cap S= \emptyset$.
We say that $N$ is an \textit{$S$-second submodule of $M$} if satisfies the equivalent conditions
of Theorem \ref{t2.333}. By an \textit{$S$-second module}, we mean
a module which is an $S$-second submodule of itself.
\end{defn}

Let $S$ be a m.c.s. of $R$. Recall that the saturation $S^*$ of $S$ is defined as $S^*=\{x \in R : x/1\  is \ a\ unit  \ of\ S^{-1}R \}$. It is obvious that $S^*$ is a m.c.s. of $R$ containing $S$ \cite{Gr92}.

\begin{prop}\label{p2.2}\cite[Proposition 2.5]{MR4401391}
Let $S$ be a m.c.s. of $R$ and $M$ be an R-module. Then we have the following.
\begin{itemize}
\item [(a)] If $N$ is a second submodule of $M$ such that $S \cap Ann_R(N)=\emptyset$, then $N$ is an $S$-second submodule of $M$. In fact if $S \subseteq u(R)$ and $N$ is an $S$-second submodule of $M$, then $N$ is a second submodule of $M$.
\item [(b)] If $S_1 \subseteq S_2$ are m.c.s.s of $R$ and $N$ is an $S_1$-second submodule of $M$, then $N$ is an $S_2$-second submodule of $M$ in case $Ann_R(N) \cap S_2=\emptyset$.
\item [(c)] $N$ is an $S$-second submodule of $M$ if and only if $N$ is an $S^*$-second submodule of $M$
\item [(d)] If $N$ is a finitely generated $S$-second submodule of $M$, then $S^{-1}N$ is a second submodule of $S^{-1}M$
\end{itemize}
\end{prop}

\begin{cor}\label{c2.2}\cite[Corollary 2.6]{MR4401391}
Let $M$ be an R-module and set $S=\{1\}$. Then  every second submodule of $M$ is an $S$-second submodule of $M$.
\end{cor}

The following examples show that the converses of Proposition \ref{p2.2} (a) and (d) are not true in general.
\begin{ex}\label{e2.2}\cite[Example 2.7]{MR4401391}
Take the $\Bbb Z$-module $M=\Bbb Z_{p^\infty}\oplus \Bbb Z_2$ for a prime number $p$. Then $2(\Bbb Z_{p^\infty}\oplus \Bbb Z_2)=\Bbb Z_{p^\infty}\oplus 0$ implies that $M$ is not a second $\Bbb Z$-module.  Now,
take the m.c.s. $S =\Bbb Z\setminus \{0\}$ and put $s=2$. Then $2rM=\Bbb Z_{p^\infty}\oplus 0=2M$ for all $r \in \Bbb Z$ and so  $M$ is an $S$-second $\Bbb Z$-module.
\end{ex}

\begin{ex}\label{e22.2}\cite[Example 2.8]{MR4401391}
Consider the $\Bbb Z$-module $M=\Bbb Q \oplus \Bbb Q$, where $\Bbb Q$ is the field of rational numbers. Take the submodule
$N = \Bbb Z \oplus 0$ and the m.c.s.  $S =\Bbb Z\setminus \{0\}$. Then one can see that $N$ is not an $S$-second submodule of $M$. Since $S^{-1}\Bbb Z =\Bbb Q$
is a field, $S^{-1}(\Bbb Q \oplus \Bbb Q)$ is a vector space so that a non-zero submodule $S^{-1}N$ is a second submodule of $S^{-1}(\Bbb Q \oplus \Bbb Q)$.
\end{ex}

\begin{prop}\label{p2.4}\cite[Proposition 2.9]{MR4401391}
Let $M$ be an $R$-module and $S$ be a m.c.s. of $R$. Then the following statements hold.
\begin{itemize}
\item [(a)] If $N$ is an $S$-second submodule of $M$, then $Ann_R(N)$ is an $S$-prime ideal of $R$.
\item [(b)] If $M$ is a comultiplication $R$-module and $Ann_R(N)$ is an $S$-prime ideal of $R$, then $N$ is an $S$-second submodule of $M$.
\end{itemize}
\end{prop}

An $R$-module $M$ satisfies the \emph{double annihilator
conditions} (DAC for short)  if for each ideal $I$ of $R$
we have $I=Ann_R(0:_MI)$ \cite{Fa95}. An
$R$-module $M$ is said to be a \emph{strong comultiplication module} if $M$ is
a comultiplication $R$-module and satisfies the DAC conditions \cite{AF09}.

\begin{thm}\label{t2.5}\cite[Theorem 2.10]{MR4401391}
Let $M$ be a strong comultiplication $R$-module and $N$ be a submodule of $M$ such that $Ann_R(N) \cap S=\emptyset$, where $S$ is a m.c.s. of $R$. Then the following are equivalent:
\begin{itemize}
\item [(a)] $N$ is an $S$-second submodule of $M$;
\item [(b)] $Ann_R(N)$ is an $S$-prime ideal of $R$;
\item [(c)] $N=(0:_MI)$ for some $S$-prime ideal $I$ of $R$ with $Ann_R(N) \subseteq I$.
\end{itemize}
\end{thm}

Let $R_i$ be a commutative ring with identity, $M_i$ be an $R_i$-module for each $i = 1, 2,..., n$, and $n \in \Bbb N$. Assume that
$M = M_1\times M_2\times ...\times M_n$ and $R = R_1\times R_2\times ...\times R_n$. Then $M$ is clearly
an $R$-module with componentwise addition and scalar multiplication. Also,
if $S_i$ is a multiplicatively closed subset of $R_i$ for each $i = 1, 2,...,n$,  then
$S = S_1\times S_2\times ...\times S_n$ is a multiplicatively closed subset of $R$. Furthermore,
each submodule $N$ of $M$ is of the form $N = N_1\times N_2\times...\times N_n$, where $N_i$ is a
submodule of $M_i$.
\begin{thm}\label{t2.6}\cite[Theorem 2.11]{MR4401391}
Let  $M = M_1 \times M_2$ be an $R = R_1 \times R_2$-module and $S = S_1\times S_2$ be a  m.c.s. of
$R$, where $M_i$ is an $R_i$-module and $S_i$ is a  m.c.s. of $R_i$ for each $i = 1, 2$. Let $N=N_1\times N_2$ be a submodule of $M$. Then the following are equivalent:
\begin{itemize}
\item [(a)] $N$ is an $S$-second submodule of $M$;
\item [(b)] $N_1$ is an $S_1$-second submodule of $M_1$ and $Ann_{R_2}(N_2) \cap S_2\not=\emptyset$ or $N_2$ is an $S_2$-second submodule of $M_2$ and $Ann_{R_1}(N_1) \cap S_1\not=\emptyset$.
\end{itemize}
\end{thm}

\begin{thm}\label{t2.7}\cite[Theorem 2.12]{MR4401391}
Let  $M = M_1 \times M_2\times ... \times M_n$ be an $R = R_1 \times R_2\times ...\times R_n$-module and $S = S_1\times S_2\times ... \times S_n$ be a  m.c.s. of
$R$, where $M_i$ is an $R_i$-module and $S_i$ is a m.c.s. of $R_i$ for each $i = 1, 2,...,n$. Let $N=N_1\times N_2\times ... \times N_n$ be a submodule of $M$. Then the following are equivalent:
\begin{itemize}
\item [(a)] $N$ is an $S$-second submodule of $M$;
\item [(b)] $N_i$ is an $S_i$-second submodule of $M_i$ for some $i \in \{1,2,...,n\}$ and $Ann_{R_j}(N_j) \cap S_j\not=\emptyset$ for all $j \in \{1,2,...,n\}-\{ i \}$.
\end{itemize}
\end{thm}

\begin{lem}\label{t22.8}\cite[Lemma 2.13]{MR4401391}
Let $S$ be a m.c.s. of $R$ and $N$ be an $S$-second submodule of an $R$-module $M$. Then the following statements hold for some $s \in S$.
\begin{itemize}
\item [(a)] $sN \subseteq \acute{s}N$ for all $\acute{s}\in S$.
\item [(b)] $(Ann_R(N):_R\acute{s}) \subseteq (Ann_R(N):_Rs)$ for all $\acute{s} \in S$.
\end{itemize}
\end{lem}

\begin{prop}\label{t2.8}\cite[Proposition 2.14]{MR4401391}
Let $S$ be a m.c.s. of $R$ and $N$ be a finitely generated submodule of $M$ such that $Ann_R(N) \cap S=\emptyset$. Then the following are equivalent:
\begin{itemize}
\item [(a)] $N$ is an $S$-second submodule of $M$;
\item [(b)] $S^{-1}N$ is a second submodule of $S^{-1}M$ and there is an $s \in S$ satisfying $sN \subseteq \acute{s}N$ for all $\acute{s}\in S$.
\end{itemize}
\end{prop}

\begin{thm}\label{t2.8}\cite[Theorem 2.15]{MR4401391}
Let $S$ be a m.c.s. of $R$ and $N$ be a submodule of an $R$-module $M$ such that $Ann_R(N) \cap S=\emptyset$. Then $N$ is an $S$-second submodule of $M$ if and only if $sN$ is a second submodule of $M$ for some $s \in S$.
\end{thm}

The set of all maximal ideals of $R$ is denoted by $Max(R)$.
\begin{thm}\label{t2.9}\cite[Theorem 2.16]{MR4401391}
Let $S$ be a m.c.s. of $R$ and $N$ be a submodule of an $R$-module $M$ such that $Ann_R(N) \subseteq Jac(R)$, where $Jac(R)$ is the Jacobson radical of $R$. Then the following statements are equivalent:
\begin{itemize}
\item [(a)] $N$ is a second submodule of $M$;
\item [(b)] $Ann_R(N)$ is a prime ideal of $R$ and $N$ is an  $(R\setminus \mathfrak{M})$-second submodule of $M$ for each $\mathfrak{M} \in Max(R)$.
\end{itemize}
\end{thm}

Now we determine all second submodules of a module over a quasilocal ring in terms of $S$-second submodules.
\begin{cor}\label{cc2.9}\cite[Corollary 2.17]{MR4401391}
Let $S$ be a m.c.s. of  a quasilocal ring $(R, \mathfrak{M})$ and $N$ be a submodule of an $R$-module $M$. Then the following statements are equivalent:
\begin{itemize}
\item [(a)] $N$ is a second submodule of $M$;
\item [(b)] $Ann_R(N)$ is a prime ideal of $R$ and $N$ is an  $(R\setminus \mathfrak{M})$-second submodule of $M$.
\end{itemize}
\end{cor}

\begin{prop}\label{t2.16}\cite[Proposition 2.18]{MR4401391}
Let $S$ be a m.c.s. of $R$ and $f : M \rightarrow \acute{M}$ be a monomorphism of R-modules. Then we have the following.
\begin{itemize}
  \item [(a)] If $N$ is an $S$-second submodule of $M$, then $f(N)$ is an $S$-second submodule of $\acute{M}$.
  \item [(b)] If $\acute{N}$ is an $S$-second submodule of $\acute{M}$ and $\acute{N} \subseteq f(M)$, then $f^{-1}(\acute{N})$ is an $S$-second submodule of $M$.
 \end{itemize}
\end{prop}

\begin{prop}\label{p2.17}\cite[Proposition 2.19]{MR4401391}
Let $S$ be a m.c.s. of $R$, $M$ a comultiplication $R$-module, and  let $N$ be an $S$-second submodule of $M$. Suppose that $N \subseteq K+H$ for some
submodules $K, H$ of $M$. Then $sN \subseteq K$ or $sN \subseteq H$ for some $s \in S$.
\end{prop}

Let $M$ be an $R$-module. The idealization $R(+)M =\{(a,m): a \in R, m \in  M\}$ of $M$ is
a commutative ring whose addition is componentwise and whose multiplication is defined as $(a,m)(b,\acute{m}) =
(ab, a\acute{m} + bm)$ for each $a, b \in R$, $m, \acute{m}\in M$ \cite{Na62}. If $S$ is a m.c.s. of $R$ and $N$ is a submodule of $M$, then $S(+)N = \{(s, n): s \in S, n \in N\}$ is a m.c.s. of $R(+)M$ \cite{DD02}.

\begin{prop}\label{p2.18}\cite[Proposition 2.20]{MR4401391}
Let $M$ be an $R$-module and let $I$ be an ideal of $R$ such that $I \subseteq Ann_R(M)$.  Then the following are equivalent:
\begin{itemize}
\item [(a)] $I$ is a second ideal of $R$;
\item [(b)] $I(+)0$ is a second ideal of $R(+)M$.
\end{itemize}
\end{prop}

\begin{thm}\label{t2.19}\cite[Theorem 2.21]{MR4401391}
Let $S$ be a m.c.s. of $R$, $M$ be an $R$-module, and $I$ be an ideal of $R$ such that $I \subseteq Ann_R(M)$ and $I \cap S=\emptyset$. Then the following are equivalent:
\begin{itemize}
\item [(a)] $I$ is an $S$-second ideal of $R$;
\item [(b)] $I(+)0$ is an $S(+)0$-second ideal of $R(+)M$;
\item [(c)] $I(+)0$ is an $S(+)M$-second ideal of $R(+)M$.
\end{itemize}
\end{thm}

One can see that if $M$ is a cotorsion-free $R$-module, then $R$ is an integral domain and $M$ is a faithful $R$-module.
In Proposition \ref{p3.12} (e), it is shown that if $M$ is a comultiplication $R$-module the reverse is true. The following example shows that sometimes the reverse of this statement may not be true.

\begin{ex}\label{e2.20}\cite[Example 2.22]{MR4401391}
 Consider the $\Bbb Z$-module $M=\prod^{\infty}_{i=1}\Bbb Z_{p^i}$, where
$p$ is a prime number. Then it is easy to see that $M$ is a faithful $\Bbb Z$-module.
But the $\Bbb Z$-module $M$ is not second since $(\bar{1},\bar{0}, \bar{0},...) \not \in pM$ and so $M \not = pM$. Therefore, $I^M_0(M)\not=M$ and so the $\Bbb Z$-module $M$ is not a
cotorsion-free module.
\end{ex}

\begin{defn}\label{d2.20}\cite[Definition 2.23]{MR4401391}
Let $M$ be an $R$-module and $S$ be a m.c.s. of $R$  with $Ann_R(M) \cap S=\emptyset$. We say that $M$ is an \emph{$S$-cotorsion-free module} in the case that we can find $s\in S$ such that  if $rM \subseteq L$, where $r \in R$ and $L$ is a completely irreducible submodule of $M$, then $sM \subseteq L$ or $rs=0$.
\end{defn}

\begin{prop}\label{p22.3}\cite[Proposition 2.24]{MR4401391}
Let $M$ be an $R$-module and $S$ be a m.c.s. of $R$.
Then the following statements are equivalent.
\begin{itemize}
  \item [(a)] $M$ is an $S$-second $R$-module.
  \item [(b)] $\mathfrak{p}=Ann_R(M)$ is an $S$-prime ideal of $R$ and the $R/\mathfrak{p}$-module $M$ is
an $S$-cotorsion-free module.
\end{itemize}
\end{prop}

\begin{thm}\label{t2.21}\cite[Theorem 2.25]{MR4401391}
Let $M$ be a module over an integral domain $R$. Then the following are equivalent:
\begin{itemize}
\item [(a)]
 $M$ is a cotorsion-free $R$-module;
\item [(b)]
$M$ is an $(R \setminus \mathfrak{p})$-cotorsion-free for each prime ideal $\mathfrak{p}$ of $R$;
\item [(c)]
$M$ is an $(R \setminus \mathfrak{M})$-cotorsion-free for each maximal ideal $\mathfrak{M}$ of $R$.
\end{itemize}
\end{thm}

\begin{thm}\label{t2.10}\cite[Theorem 2.26]{MR4401391}
Let $S$ be a m.c.s. of $R$ and $M$ be a finitely generated comultiplication $R$-module with $Ann_R(M) \cap S=\emptyset$. Then the following statements are equivalent:
\begin{itemize}
\item [(a)] Each non-zero submodule of $M$ is $S$-second;
\item [(b)] $M$ is a simple $R$-module.
\end{itemize}
\end{thm}

\begin{cor}\label{c2.110}\cite[Corollary 2.27]{MR4401391}
Let $S$ be a m.c.s. of $R$. If $M$ is a finitely generated multiplication and comultiplication $R$-module with $Ann_R(M) \cap S=\emptyset$, then the following statements are equivalent:
\begin{itemize}
\item [(a)] Each non-zero submodule of $M$ is $S$-second;
\item [(b)] $M$ is a simple $R$-module;
\item [(c)] Each proper submodule of $M$ is an $S$-prime submodule of $M$.
\end{itemize}
\end{cor}

\begin{ex}\label{e2.10}\cite[Example 2.28]{MR4401391}
Consider the $\Bbb Z$-module $\Bbb Z_n$. Take $S=\Bbb Z -0$. We know that $\Bbb Z_n$ is a finitely generated multiplication and comultiplication $\Bbb Z$-module. Then by Corollary \ref{c2.110}, if $n$ is not a prime number, the $\Bbb Z$-module $\Bbb Z_n$ has a non-zero submodule which is not $S$-second and a proper submodule which is not $S$-prime.
\end{ex}

\begin{thm}\label{e2.10}\cite[Theorem 3.6]{MR4453905}
Every $S$-comultiplication $S$-torsion-free module is an $S$-cyclic module
\end{thm}

\begin{thm}\label{e2.10}\cite[Theorem 4.5]{MR4453905}
 Let $N$ be a submodule of an $R$-module $M$ with $Ann_R(N) \cap  S = \emptyset$. The following assertions are equivalent.
 \begin{itemize}
\item [(a)] $N$ is an $S$-second submodule.
\item [(b)] There exists $s \in  S$ such that for each $a \in R$, the homothety $ a:N\hookrightarrow N$ is either $S$-zero or $S$-surjective
with respect to $s \in  S$.
\item [(c)] There exists a fixed $s \in S$ so that for each $a \in  R$, either $saN = 0$ or $sN \subseteq aN$.
\end{itemize}
\end{thm}

\begin{thm}\label{e2.10}\cite[Theorem 4.6]{MR4453905}
 Let $M$ be an $S$-comultiplication module. The following statements are equivalent.
 \begin{itemize}
\item [(a)] $N$ is an $S$-second submodule of $M$.
\item [(b)] $Ann_R(N)$ is an $S$-prime ideal of $R$ and there exists $s \in S$ such that $sN \subseteq \acute{s}N$ for each $\acute{s} \in S$.
\end{itemize}
\end{thm}

\begin{thm}\label{e2.10}\cite[Theorem 4.8]{MR4453905}
Let $M$ be an $S$-comultiplication $R$-module and let $N$ be an $S$-second submodule of $M$. If
$N \subseteq N_1 + N_2 +\ldots + N_m$ for some submodules $N_1,N_2,\ldots ,N_m$ of $M$, then there exists $s \in S$ such that
$sN \subseteq N_i$ for some $1\leq i \leq m$.
\end{thm}
\section{Graded second submodule}
\noindent
\begin{defn}\cite[Definition 3.14]{MR2832354}
Let $R$ be a $G$-graded ring, $M$ be a graded $R$-module.  A non zero graded submodule $N$ of $M$ is said to be a \textit{graded second submodule} of $M$ if $rN=0$ or $rN=N$ for every $r\in~h(R)$.
\end{defn}

\begin{rem}
\cite[Remark 2.1]{AF1089891} It is clear that every
second $R$-module which is a graded module is a gr-second
$R$-module but the converse is not true in general. For example
if we take $R=K[x,  x^{-1}](=K[x]_x)$, where $K$ is a field and $x$ is an
indeterminate, graded in the obvious way, $R$ as an $R$-module
is graded simple (see \cite[1.5.14(c)]{BH96}). Hence $R$ is a gr-second
$R$-module. But $R$ is not a secondary $R$-module by \cite{SH86}. Hence $R$ is not a
second $R$-module.
 \end{rem}

\begin{prop}
\cite[Proposition 2.3]{AF1089891} Let $M$ be a graded $R$-module. Then the
following hold.
\begin{itemize}
  \item [(a)] If $S$ is a gr-secondary submodule of $M$,
  then $S$ is a gr-second if and only if $Ann_R(S)$ is a gr-prime ideal of $R$.
  \item [(b)] Let $K$ be a graded submodule of a $P$-gr-second
  module $M$. Then $K$ is a $P$-gr-secondary submodule if and only if
  $K$ is a $P$-gr-second submodule.
  \item [(c)] If $S$ is a gr-minimal submodule of $M$, then $S$ is a
  gr-second submodule of $M$.
\end{itemize}
 \end{prop}

\begin{prop}
\cite[Proposition 2.4]{AF1089891} Let $R$ be a graded ring and $P$ be a
gr-prime ideal of $R$. Then we have the following.
\begin{itemize}
  \item [(a)] The sum of $P$-gr-second $R$-modules  is a $P$-gr-second
  module.
  \item [(b)] The product of $P$-gr-second $R$-modules is a $P$-gr-second
  module.
  \item [(c)]  Every non-zero gr-quotient of a $P$-gr-second $R$-module is a $P$-gr-second
  module.
\end{itemize}
 \end{prop}

\begin{lem}
\cite[Lemma 2.5]{AF1089891} Let $P$ be a graded
prime ideal of $R$ and let $S$ be a graded non-zero submodule of
a graded $R$-module $M$. Then the following are equivalent.
\begin{itemize}
  \item [(a)] $S$ is a $P$-gr-second submodule of $M$.
  \item [(b)] $Ann_R(S)=W^{gr}(S)=P$, where
   $$
  W^{gr}(M)= \{ a \in h(R): \\\ the \\\ homothety \\\
  M \stackrel {a} \rightarrow M \\\ is  \\\ not \\\ surjective \}.
  $$
\end{itemize}
 \end{lem}

\begin{thm}
\cite[Theorem 2.6]{AF1089891} Let $M$ be a graded $R$-module and let $S$ be
a non-zero graded submodule of $M$ with $Ann_R(S)=P$ is a graded
prime ideal of $R$. Then the following are equivalent.
\begin{itemize}
  \item [(a)] $S$ is a $P$-gr-second submodule of $M$.
  \item [(b)] $S$ is a gr-divisible $R/P$-module.
  \item [(c)] $rS=S$ for all $r \in h(R)-P$.
  \item [(d)] $IS=S$ for all graded ideals $I$ with $I \subseteq P$.
  \item [(e)] $W^{gr}(S)=P$.
\end{itemize}
 \end{thm}

\begin{defn}
\cite[Definition 2.1]{AF1089891}We shall call a graded submodule $N$ of $M$
a  \emph{minimal $P$-gr-secondary} (resp.  \emph{$P$-gr-second})
submodule of $M$. If $N$ is a $P$-gr-secondary (resp.
$P$-gr-second) submodule which is not strictly contains any other
$P$-gr-secondary
(resp. $P$-gr-second) submodule of $M$.
 \end{defn}

\begin{thm}
\cite[Theorem 2.7]{AF1089891} Let $M$ be a graded $R$-module. Then the
submodule $N$ of $M$ is minimal $P$-gr-secondary
if and only if $N$ is a minimal $P$-gr-second submodule of $M$.
 \end{thm}

\begin{thm}
\cite[Theorem 2.8]{AF1089891} Let $M$ be a gr-prime module. Then the
following are equivalent.
\begin{itemize}
  \item [(a)] $M$ is a gr-second module.
  \item [(b)] $M$ is an gr-injective $R/Ann_R(M)$-module.
\end{itemize}
 \end{thm}

\begin{prop}
\cite[Proposition 2.9]{AF1089891} Let $M$ be a graded $R$-module and let $N$ be a graded submodule
 of $M$ . Then
 \begin{itemize}
   \item [(a)] If $M$ is a gr-primary module and $N$ is a
    gr-second submodule of $M$, then $N$ is $Ann_R(N)$-gr-primary.
   \item [(b)] If $M$ is a gr-prime module and $N$ is a
   gr-second submodule of $M$, then $rN=rM \cap N$ for each $r \in h(R)$.
   \item [(c)] If $Ann_R(N)$ is a gr-prime ideal of $R$
   and $N$ is a gr-minimal in the set of all graded
   submodules $K$ of $M$ such that $Ann_R(K)=Ann_R(N)$,
   then $N$ is a graded second submodule of $M$.
 \end{itemize}
 \end{prop}

\begin{thm}
\cite[Theorem 2.10]{AF1089891} Let $E$ be a graded injective cogenerator of
$R$ and let $N$ be a graded submodule of a graded $R$-module $M$.
Then $N$ is a gr-prime submodule of $M$ if and only if
$HOM_R(M/N,E)$ is
a gr-second $R$-module.
 \end{thm}

\begin{thm}
\cite[Theorem 2.11]{AF1089891} Let $R$ be a graded integral domain which is not a
gr-field and $K$ the gr-field of quotient of $R$. Then the $R$-module $K$
has no gr-minimal submodule and $K$ is the only gr-second submodule
of $K$.
 \end{thm}

\begin{prop}
\cite[Proposition 3.15]{MR2832354} Let $M$ be a graded $R$-module and let $N$
be a graded submodule of $M$. Then we have the following.
\begin{itemize}
  \item [(a)] If $N$ is a gr-second submodule of $M$, then
  $Ann_R(N)$ is a gr-prime ideal, $P$ say, of $R$, we say that $N$
  is $P$-gr-second submodule of $M$.
  \item [(b)] If $M$ is a gr-comultiplication $R$-module and
  $Ann_R(N)$ is a gr-prime ideal of $R$, then $N$ is a gr-second
  submodule of $M$.
\end{itemize}
 \end{prop}

\begin{thm}
\cite[Theorem 3.16]{MR2832354} Let $M$ be a Noetherian gr-comultiplication
$R$-module. Then we have the following.
\begin{itemize}
  \item [(a)] $M$ has a finite number of gr-second submodules.
  \item [(b)] Every gr-second submodule of $M$ is a gr-minimal submodule of $M$.
\end{itemize}
 \end{thm}

\section{$I$-second submodules of a modules}
\begin{thm}\label{t1.991}\cite[Theorem 2.3]{MR4049591}
 Let $I$ be an ideal of $R$. For a non-zero submodule $S$ of an $R$-module $M$ the following statements are equivalent:
\begin{itemize}
  \item [(a)] For each $r \in R$, a submodule $K$ of $M$, $r \in (K:_RS) \setminus (K:_R(S:_MI))$ implies that $S \subseteq K$ or $r \in Ann_R(S)$;
  \item [(b)] For each $r \not \in (rS:_R(S:_MI))$, we have $rS=S$ or $rS=0$;
  \item [(c)] $(K:_RS) = Ann_R(S) \cup (K:_R(S:_MI))$, for any submodule $K$ of $M$ with $S \not \subseteq K$;
  \item [(d)] $(K:_RS) = Ann_R(S)$ or $(K:_RS)=(K:_R(S:_MI))$, for any submodule $K$ of $M$ with $S \not \subseteq K$.
 \end{itemize}
\end{thm}

\begin{defn}\label{d1.1}\cite[Definition 2.4]{MR4049591}
Let $I$ be an ideal of $R$. We say that a non-zero submodule $S$ of an $R$-module $M$ is an $I$-\emph{second submodule} of $M$ if satisfies the equivalent conditions of Theorem \ref{t1.991}. This can be regarded as a dual notion of the $I$-prime submodule. In case, $I=0$ we say that $S$ is a \emph{weak second submodule} of $M$.
\end{defn}

Let $I$ be an ideal of $R$. Clearly every second submodule is an $I$-second submodule. But the converse is not true in general as we see in the following example.
\begin{ex}\label{e222.14}\cite[Example 2.5.]{MR4049591}
\begin{itemize}
\item [(a)]
If $I=0$, then every module is an $I$-second submodule of itself but every module is not a second module. For example, the $\Bbb Z$-module $\Bbb Z$ is weak second which is not second.
\item [(b)]
Consider the $\Bbb Z$-module $\Bbb Z_{12}$. Take $I = 4\Bbb Z$ as an ideal of $\Bbb Z$ and $S = \bar{3}\Bbb Z_{12}$ as a submodule of $\Bbb Z_{12}$. Then $S$ is an $I$-second submodule of $\Bbb Z_{12}$. But  $S$ is not second submodule.
\end{itemize}
\end{ex}

\begin{ex}\label{e2.2}\cite[Example 2.6.]{MR4049591}
Let $I$ be an ideal of $R$ and $S$ a non-zero submodule of an $R$-module $M$. If for each $r \in R$, a completely irreducible submodule $L$ of $M$, $r \in (L:_RS) \setminus (L:_R(S:_MI))$ implies that $S \subseteq L$ or $r \in Ann_R(S)$ we can not conclude that $S$ is an $I$-second submodule of $M$. For example, consider $\Bbb Z$ as a $\Bbb Z$-module. Then $2\Bbb Z$ satisfies the mention condition above but it is not  an $I$-second submodule of $\Bbb Z$ for ideal $I=4 \Bbb Z$ of $\Bbb Z$.
\end{ex}

Let $I$ be an ideal of $R$ and $M$ be an $R$-module. If $I=R$, then every submodule is an $I$-second submodule. So in the rest of this paper we can assume that $I \not =R$.

\begin{thm}\label{t1.3}\cite[Theorem 2.7]{MR4049591}
 Let $M$ be an $R$-module. Then we have the following.
  \begin{itemize}
    \item [(a)] Let $I, J$ be ideals of $R$ such that $I \subseteq J$. If $S$ is an $I$-second submodule of $M$, then $S$ is  an $J$-second submodule of $M$. In particular, every weak second submodule is an $I$-second submodule for each ideal $I$ of $R$.
    \item [(b)] If $S$ an $I$-second submodule of $M$ which is not second, then
  $Ann_R(S)(S:_MI) \subseteq S$.
\end{itemize}
\end{thm}

\begin{thm}\label{t1.5}\cite[Theorem 2.8]{MR4049591}
Let $I$ be an ideal of $R$, $M$ an $R$-module, and $S$ be a submodule of $M$. Then we have the following.
\begin{itemize}
  \item [(a)] If $S$ is an $I$-second submodule of $M$ such that $Ann_R((S:_MI))\subseteq IAnn_R(S)$, then $Ann_R(S)$ is an $I$-prime ideal of $R$.
  \item [(b)] If $M$ is a comultiplication $R$-module and $Ann_R(S)$ is an $I$-prime ideal of $R$, then $S$ is an $I$-second submodule of $M$.
\end{itemize}
\end{thm}

\begin{cor}\label{c1.5}\cite[Corollary 2.9]{MR4049591}
Let $M$ an $R$-module and $S$ be a submodule of $M$. Then we have the following.
\begin{itemize}
  \item [(a)] If $M$ is faithful and $S$ is a weak second submodule of $M$, then $Ann_R(S)$ is a weakly prime ideal of $R$.
  \item [(b)] If $M$ is a comultiplication $R$-module and $Ann_R(S)$ is a weakly prime ideal of $R$, then $S$ is a weak second submodule of $M$.
\end{itemize}
\end{cor}

The following example shows that the condition ``$M$ is a comultiplication $R$-module" in  Corollary \ref{c1.5} (b) can not be omitted.
\begin{ex}\label{e22.14}\cite[Example 2.10]{MR4049591}
Let $R= \Bbb Z$, $M=\Bbb Z\oplus \Bbb Z$, and $S=2\Bbb Z \oplus 0$. Then $M$ is not a comultiplication $R$-module. Clearly,  $Ann_R(S)=0$ is a  weakly prime ideal of $R$. But $S$ is not a weak second submodule of $M$.
\end{ex}
\begin{prop}\label{t1.5}\cite[Proposition 2.11]{MR4049591}
Let $I$ be an ideal of $R$ and $M$ be an $R$-module. Let $N$ be an $I$-second
submodule of $M$. Then we have the following statements.
\begin{itemize}
  \item [(a)] If $K$ is a submodule of $M$ with $K \subset N$, then $N/K$ is an $I$-second submodule of $M/K$.
  \item [(b)] Let $N$ be a finitely generated submodule of $M$ and $S$ be a  multiplicatively closed subset of $R$ with $Ann_R(N) \cap S=\emptyset$. Then $S^{-1}N$ is an $S^{-1}I$-second submodule of $S^{-1}M$.
\end{itemize}
\end{prop}

\begin{thm}\label{t1.6}\cite[Theorem 2.12]{MR4049591}
Let $M$ be a primary $R$-module. Then every proper weak second submodule of $M$ is a primary submodule of $M$.
\end{thm}

\begin{prop}\label{p1.6}\cite[Proposition 2.13]{MR4049591}
Let $I$ be an ideal of $R$, $M$ and $\acute{M}$ be $R$-modules, and let $f : M\rightarrow \acute{M}$ be an $R$-monomorphism. If $\acute{N}$ is an $I$-second submodule of $\acute{M}$ such that $\acute{N} \subseteq Im(f)$, then $f^{-1}(\acute{N})$ is an $I$-second submodule of $M$.
\end{prop}

\begin{lem}\label{l1.7}\cite[Lemma 2.14]{MR4049591}
Let $R = R_1 \times R_2$ be a decomposable ring, $I=I_1 \times I_2$ an ideal of $R$, and $M = M_1 \times M_2$ be an $R$-module, where $M_1$ is an $R_1$-module and $M_2$ is an $R_2$-module. If $(0:_{M_2}I_2)\not=0$ and  $S_1$ is a non-zero $R_1$-submodule of $M_1$, then the following statements are equivalent:
\begin{itemize}
  \item [(a)] $S_1$ is a second $R_1$-submodule of $M_1$;
  \item [(b)] $S_1 \times 0$ is a second $R$-submodule of $M = M_1 \times M_2$;
  \item [(c)] $S_1 \times 0$ is an $I$-second $R$-submodule of $M = M_1 \times M_2$.
\end{itemize}
\end{lem}

\begin{thm}\label{p1.8}\cite[Theorem 2.15]{MR4049591}
Let $R = R_1 \times R_2$ be a decomposable ring and $M = M_1 \times M_2$
be an $R$-module, where $M_1$ is an $R_1$-module and $M_2$ is an $R_2$-module. Let $I$ be an ideal of $R$ such that $(0:_{M_1}I_1) \not =0$ and $(0:_{M_2}I_2)\not=0$. If $S =S_1\times S_2$ is an $I$-second $R$-submodule of $M = M_1\times M_2$, then either $(S:_MI)=S$ or $S$ is a second submodule of $M$.
\end{thm}

\begin{ex}\label{e2.14}\cite[Example 2.16]{MR4049591}
Let $R_1 = R_2 = M_1 = M_2 =S_1= \Bbb Z_6$. Then by Theorem \ref{p1.8},  $S_1 \times 0$ is not a weak second submodule of $ M_1 \times M_2$.
\end{ex}

\begin{thm}\label{t1.9}\cite[Theorem 2.17]{MR4049591}
Let $I$ be an ideal of $R$, $M_1$, $M_2$ be $R$-modules, and let $N$ be a submodule of $M_1$. Then $N \oplus 0$ is an $I$-second submodule of $M_1\oplus M_2$ if and only if $N$ is an $I$-second submodule of $M_1$ and for $r\in (rN:_R(N:_{M_1}I))$, $rN \not =0$, and $rN \not =N$, we have $r\in Ann_R((0:_{M_2}I))$.
\end{thm}

\begin{cor}\label{c1.10}\cite[Corollary 2.18]{MR4049591}
Let $I$ and $\mathfrak{p}$ be  ideals of $R$, $M_1$, $M_2$ be $R$-modules, and let $N$ be a submodule of $M_1$. Let $S_i$ ($1\leq i \leq n$) be $\mathfrak{p}$-secondary submodules of $M_1$ with $\sum^n_{i=1}S_i=(N:_{M_1}I)$. If $N$ is an $I$-second submodule of $M_1$ and $\mathfrak{p}\subseteq Ann_R((0:_{M_2}I))$, then $N \oplus 0$ is an $I$-second submodule of $M_1 \oplus M_2$.
\end{cor}

\begin{thm}\label{t1.11}\cite[Theorem 2.19]{MR4049591}
Let $I$ be an ideal of $R$ and $M$ be an $R$-module. Then we have the following.
\begin{itemize}
  \item [(a)] If $\cap^\infty_{n=1}I^nM=0$ and every proper submodule of $M$ is $I$-prime, then every non-zero submodule of $M$ is $I$-second.
  \item [(b)] If $\sum^\infty_{n=1}(0:_MI^n)=M$ and every non-zero submodule of $M$ is $I$-second, then every proper submodule of $M$ is $I$-prime.
\end{itemize}
\end{thm}

\begin{cor}\label{c1.12}\cite[Corollary 2.20]{MR4049591}
Let $I$ be an ideal of $R$ and $M$ be an $R$-module. Then  every proper submodule of $M$ is weakly prime if and only if every non-zero submodule of $M$ is weak second.
\end{cor}

\begin{cor}\label{c1.12}\cite[Corollary 2.21]{MR4049591}
Let $(R, m)$ be a local ring and $M$ be an $R$-module. Then we have the following.
 \begin{itemize}
  \item [(a)] If $M$ is a Noetherian $R$-module and every proper submodule of $M$ is $I$-prime, then every non-zero submodule of $M$ is $I$-second.
  \item [(b)] If $M$ is an Artinian $R$-module and every non-zero submodule of $M$ is $I$-second, then every proper submodule of $M$ is $I$-prime.
\end{itemize}
\end{cor}
\section{Fuzzy second submodules}
In this section, first we recall some basic definitions and remarks which are needed in the sequel.

A fuzzy subset $\mu$ of a non-empty set $X$ is defined as a map from $X$ to the unit interval $I := [0, 1]$. We denote the set of all fuzzy subsets of $X$ by $F(X)$. Let $\mu, \lambda \in F(X)$. Then the inclusion $\mu \subseteq \lambda$ (resp. $\mu \subset \lambda$) is denoted by $\mu(x)\leq \lambda(x)$ (resp. $\mu(x) < \lambda(x)$) for all $x \in X$.

We write $\wedge$ and $\vee$ for infimum and supremum.

Let $\mu, \lambda \in F(X)$. Then $\mu \cap \lambda$ and $\mu \cup \lambda$ are defined as follows: \\ \begin{center} $(\mu\cap\lambda)(x)= \mu(x) \wedge \lambda(x)$,\\ $(\mu \cup \lambda)(x) = \mu(x) \vee \lambda(x)$,  for all $x \in X$.\end{center}

Let $\mu \in F(X)$. Then $\mu$ has sup property if every subset of $\mu(x)$ has a maximal element.

Let $\mu \in F(X)$ and $t\in I$. Then the set $\mu_t = \{x\in X |\:\mu(x) \leq t\}$ is called the \emph{$t$-level subset} of $X$ with respect to $\mu$ and $\mu^*=\{x \in X, \mu(x)>\mu(0)\}$ is called \emph{the support} of $\mu$. Also if $\mu \in F(X)$, then $\mu_* = \{x \in X,\:\mu(x)=\mu(0)\}$.

\begin{defn}\label{charc} (See \cite{2}.)Let $Y \subseteq X$ and $t \in I$. Then
 \begin{center}
$t_Y(x) =$ $\left\{
   \begin{array}{ll}
     t & \hbox{if  $x \in Y$ ;} \\
     0 & \hbox{Otherwise.}
   \end{array}
 \right.
$

\end{center}

In particular, if $Y = \{y\}$, then $t_{y}$ is often refereed to a fuzzy singleton point (or fuzzy point). Moreover, $1_Y$ is refereed as the characteristic function of $Y$.

If $t_y$ is a fuzzy singleton point and $t_y \subseteq \mu \in F(X)$, we write $t_y \in \mu$.
\end{defn}

\begin{defn} ( See \cite{2}.)
\begin{itemize}
\item Let $f$ be a mapping from $X$ into $Y$ such that $\mu \in F(X)$ and $\nu \in F(Y)$. Then $f(\mu) \in F(Y)$ and $f^{-1}(\nu) \in F(X)$, defined by
$\forall y \in Y$,
\begin{center}
$f(\mu(y))=$ $\left\{
             \begin{array}{ll}
              \vee \{\mu(x)\,\mid\,x\in X,\, f(x) = y\} & \hbox{if $f^{-1}(y) \neq \emptyset$;} \\
               0, & \hbox{otherwise.}
             \end{array}
           \right.$
\end{center}

and $\forall x \in X$,\\ \begin{center}
$f^{-1}(\nu)(x) = \nu(f(x))$
\end{center}

\vspace{0.2cm}

\item A fuzzy subset $\xi$ of a ring R is
called a \emph{fuzzy ideal} of $R$ if it satisfies the following
properties:
\item [(i)] $\xi(x - y) \geq \xi(x)\wedge \xi(y))$, for all $x,y \in R$; and
\item[(ii)] $\xi(xy) \geq \xi(x) \vee \xi(y))$, for all $x, y \in R$.

\vspace{0.2cm}
\item A fuzzy subset $\mu \in F(M)$ is called a \emph{fuzzy submodule} if
\item[(i)] $\mu(\theta) = 1$,
\item[(ii)] $\mu(rx) \geq \mu(x)$, for all $r \in R$ and $x \in M$,
\item[(iii)] $\mu(x+y) \geq \mu(x) \wedge \mu(y)$, for all $x, y \in M$.

\vspace{0.2cm}

In the following, we denote the set of fuzzy submodules (resp., fuzzy ideals) of $M$ (resp,. of $R$) by $FS(M)$ (resp,. $FI(R)$. The zero fuzzy submodule of $M$(resp., fuzzy ideal of $R$) is $1_\theta$ (resp., $1_0$).

\vspace{0.2cm}

\item If $\lambda \in FS(M)$, then $\lambda_{*} = \{x \in M |\,\lambda(x) = 1\}$ is a submodule of $M$.

\vspace{0.2cm}

\item Let $\mu \in FI(R)$. Then $\Re(\mu) \in FI(R)$, defined by $\Re(\mu)(x) = \vee_{n \in \Bbb N}\mu(x^{n})$ $\forall x\in R$, is called the $\Re-radical$ of $\mu$.

\vspace{0.2cm}

\item Let $\xi \in FI(R)$ and $\mu, \nu \in FS(M)$. Define $\mu+\nu, \xi.\mu \in F(M)$ as follows \begin{center}$(\mu+\nu)(x)  = \vee \{\mu(y)\wedge\nu(z)\,|\,y,z \in M, y+z = x\}$,\\
$(\xi.\mu)(x) = \vee\{\xi(r) \wedge \mu(y)\,|\,r\in R,\, y\in M,\, ry =x\}$, for all $x \in M$.
\end{center}

\vspace{0.2cm}
\item For $\mu, \nu \in FS(M)$ and $\zeta \in LI(R)$, define $(\mu : \nu)\in FI(R)$ and $(\mu : \zeta) \in FS(M)$ as follows:
\begin{center}
$(\mu : \nu) = \cup\{\eta\,|\,\eta \in FI(R),\: \eta.\nu \subseteq \mu\}$\\
$(\mu : \zeta)= \cup\{\xi\,|\,\xi \in F(M), \:\xi.\zeta \subseteq \mu\}$
\end{center}
\vspace{0.2cm}
\item Let $\mu \in LI(R)$. Then $\mu$ is called a \emph{maximal fuzzy ideal} of $R$ if $\mu$ is a maximal element in the set of all non-constant fuzzy ideal of $R$ under ponitwise partial ordering .

\end{itemize}
\end{defn}

We recall that for every submodule $N$ of $M$, $1_N$ is the characterization function of $N$ by Definition \ref{charc}.

\begin{defn} \cite{3}
A non-constant  fuzzy submodule $\mu$ of $M$ is said to be prime fuzzy (or fuzzy prime) if for $\zeta \in FI(R)$ and $\nu \in FS(M)$ such that $\zeta.\nu \subseteq \mu$, then either $\nu \subseteq \mu$ or $\zeta \subseteq (\mu : 1_M)$. If $M = R$, then $\mu$ is  said to be a fuzzy prime ideal. In this case, we have $(\mu : 1_M) = (\mu : 1_R) = \mu$.
\end{defn}

\begin{rem} (\cite[Theorem 4.5.2]{2}.)
Let $\mu, \nu \in FS(M)$ and $\zeta \in FI(R)$. Then
\item[(i)] $(\mu : \nu) = \cup \{t_r\,|\,r\in R\, , t_r.\nu \subseteq \mu\}$;
\item[(ii)] $(\mu : \zeta) = \cup\{s_x\,|\,x\in M\, , s_x.\zeta \subseteq \mu\}$.
\end{rem}

\begin{rem}\label{a}
(\cite[Theorem 3.5]{3}.)
Let $\mu \in FS(M)$. Then $\mu$ is a fuzzy prime submodule of $M$ if and only if it satisfies the following conditions:
\item(a) $\mu_*$ is a prime submodule of $M$,
\item(b) $(\mu : 1_M)(1)$ is a prime element in I,
\item(c) $t_r.s_x \in \mu,\, r\in R,\,x\in M$, and $t, s \in I$ $\Rightarrow$ either $t_r \in (\mu : 1_M)$ or $s_x \in \mu$.
\end{rem}

\begin{rem}\label{b}(\cite[Theorem 3.6]{3}.)
Let $\mu$ be a fuzzy prime submodule of $M$. Then $(\mu : 1_M)$ is a fuzzy prime ideal.
\end{rem}

We recall that for every $r \in R$, $1_r$ is a fuzzy point of $R$ by Definition \ref{charc}, where $t =1$.

\begin{rem}\label{c}(\cite[Corollary 2.1]{1}.)
Let $1_\theta \neq \mu \in FS(M)$. Then for each $a \in R$ and each $t \in (0, 1]$, $t_a \in (1_\theta : \mu)$. That is $t_a.\mu \subseteq 1_\theta$ if and only if $1_a.\mu \subseteq 1_\theta$. Consequently, $Im(1_\theta : \mu) = \{0, 1\}$.
\end{rem}

\begin{rem} \label{d}(\cite[Remark 2.1]{1}.)
Let $\zeta \in FI(R)$ and let $\nu, \lambda \in FS(M)$ such that $\nu \subseteq \lambda$. Then
\begin{itemize}
\item[(a)]
If $\zeta.\lambda \subseteq \nu$, then $\zeta.\frac{\lambda}{\nu} \subseteq 1_\theta$.
\item[(b)]
If $\zeta.\frac{\lambda}{\nu} \subseteq 1_\theta$, then $\zeta.\lambda \subseteq 1_{\nu^*}$.
\end{itemize}
\end{rem}

\begin{rem}\label{e}(\cite[Lemma 2.2]{1}.)
Let $\nu, \lambda \in FS(M)$. Then ${(\nu : \lambda)}_* \subseteq (\nu_* :_R \lambda_*)$. In particular, ${(1_\theta : \lambda)}_* \subseteq (0 :_R  \lambda_*)$.
\end{rem}

\begin{defn}\cite[Definition 3.1]{HM2}
Let $\mu \in FS(M)$. We say that $\mu$ is a fuzzy second submodule of $M$ if for each $r\in R$, either $1_r.\mu = \mu$ or $1_r.\mu = 1_\theta$. 
\end{defn}

\begin{thm}\label{f}\cite[Theorem 3.2]{HM2}
Let $1_\theta \neq \mu$ be a fuzzy second submodule of $M$. Then $(1_\theta : \mu)$ is a fuzzy prime ideal.
\end{thm}

\begin{rem}\label{g}\cite[Remark 3.3]{HM2}
Let $1_\theta \neq \mu$ be a fuzzy second submodule. Then by Theorem \ref{f}, $\xi:= (1_\theta : \mu)$ is a fuzzy prime ideal of $R$. In this case, we say $\mu$ is $\xi$\emph{-fuzzy} \textit{second}.
\end{rem}

\begin{defn}\cite[Definition 3.4]{HM2}
Let $\mu \in FS(M)$. Then we define $W(\mu)$ as follows
\begin{center}
$W(\mu) = \cup\{1_r|\:r\in R,\:1_r.\mu \neq \mu\}$.
\end{center}
\end{defn}

\begin{prop}\label{h}\cite[Proposition 3.5]{HM2}
The following assertions are equivalent;
\begin{itemize}
\item [(a)] $\mu$ is a $\xi$-fuzzy second submodule of $M$.

\item[(b)] $(1_\theta : \mu) = W(\mu) = \xi$.
\end{itemize}
\end{prop}

\begin{prop}\label{i} \cite[Proposition 3.6]{HM2} (a) Let $\mathfrak{p}$ be prime ideal of $R$. Then $1_N$ is a $1_\mathfrak{p}$-second fuzzy submodule of $M$ if and only if $N$ is a $\mathfrak{p}$-second submodule of $M$. (We recall that, for every submodule $N$ of $M$, $1_N$ is the characterization
function of $N$ by Definition \ref{charc}).
  \item (b) Let $\mu$ be a $\xi$-fuzzy second submodule of $M$. If $\vee\{\mu(x)\,\mid\,x\not\in \mu_{*}\} < 1$, then $\mu_*$ is a $\xi_*$-second submodule of $M$.
   \end{prop}

A fuzzy submodule $\mu$ of $M$ is called a coprimary (or secondary) fuzzy submodule if for each $r \in R$, either $1_r.\mu = \mu$ or there exists $n \in \Bbb N$ such that $1_{r^n}.\mu = 1_\theta$ (or equivalently $r \in\Re(1_ \theta : \mu)$, where $\Re(\mu)(x) = \bigvee_{n \in \Bbb N}\mu(x^n)$ for all $x \in M$ \cite{1}). Clearly, every fuzzy second submodule is a fuzzy coprimary submodule but the following example shows that the converse is not true in general.

\begin{ex}\label{j}\cite[Example 3.7]{HM2}
Let $M = \Bbb Z_{p^\infty}$. Let $n>1$ and set $N := <\frac{1}{p^n} + \Bbb Z>$. Then $N$ is a coprimary $\Bbb Z$-module which is not a second submodule. Hence  $1_N$ is a fuzzy coprimary  submodule which is not a fuzzy second submodule by Proposition \ref{i} (a).
\end{ex}

\begin{thm}\label{k}\cite[Theorem 3.8]{HM2}
Let $\mu$ be a fuzzy coprimary submodule of $M$. Then $\mu$ is a fuzzy second submodule if and only if $(1_\theta : \mu)$ is a fuzzy prime ideal of $R$.
\end{thm}

\begin{cor}\label{l}\cite[Corollary 3.9]{HM2}
Let $\lambda$ be a $\xi$-fuzzy coprimary submodule of $M$ and let $\mu$ be a $\xi$-fuzzy second submodule of $M$ which contains $\lambda$. Then $\lambda$ is a $\xi$-fuzzy second submodule.
\end{cor}

We shall call a fuzzy submodule $\mu$ of $M$ a \textit{minimal $\xi$-fuzzy coprimary} (resp. \textit{$\xi$-fuzzy second}) submodule of $M$ if $\mu$ is a $\xi$-fuzzy coprimary (resp. $\xi$-fuzzy second) submodule of $M$. (We recall that $\xi = (1_\theta : \mu)$.)

\begin{prop}\label{M}\cite[Proposition 3.10.]{HM2}
Let $\mu$ be a minimal $\xi$-fuzzy second submodule of $M$. Then $\mu$ is a minimal $\xi$-fuzzy coprimary submodule of $M$.
\end{prop}

\begin{thm}\label{N}\cite[Theorem 3.11]{HM2}
Let $\mu$ and $\eta$ be two $\xi$-fuzzy submodules of $M$. Then we have the following

\item[(a)] $\mu+\eta$ is a $\xi$-fuzzy second submodule.

\item[(b)] Let $\nu\in FS(M)$ be such that $\nu \subset \mu$ and $\nu^*\subset \mu^*$. Then the fuzzy quotient submodule $\frac{\mu}{\nu}$ is $\xi$-fuzzy second.
\end{thm}

 For a multiplicative closed subset $S$ of $R$ and for a fuzzy submodule $\mu$ of $M$, we put $S(\mu):=\cap_{s\in S}(1_s.\mu)$. Clearly, $S(\mu) \in FS(M)$ and $S(\mu) \subseteq \mu$.

\begin{thm}\label{o}\cite[Theorem 3.12]{HM2}
Let $f : M\rightarrow S^{-1}M$ be the natural homomorphism $x\rightarrow x/1$ for all $x \in M$. If $\mu$ is $\xi$-fuzzy second submodule of $M$, then the following are hold.
\item[(a)] If $S\cap \xi_* \neq \emptyset$, then $S^{-1}\mu = 1_\theta$.
\item[(b)] If $S\cap \xi_* = \emptyset$, then for each $s \in S$, we have $1_s.\mu = \mu$ and so $S(\mu) = \mu$.
\item[(c)] If $\mu$ has the sup property, then either $S^{-1}\mu = 1_\theta$ or $S^{-1}\xi$-fuzzy second submodule with $(1_\theta : S^{-1}\mu) = S^{-1}\xi$.
\end{thm}

\begin{rem}\label{p}\cite[Remark 3.13]{HM2}
In \cite{1}, Theorem 3.2 says that if $f :M\rightarrow N$ is a monomorphism of non-zero submodules such that $N/Im(f)$ is a torsion free $R$-module and if $\mu$ is a $\xi$-fuzzy coprimary submodule of $M$, then so is $f(\mu)$. But the given proof is valid just for fuzzy second submodules not for fuzzy coprimary submodules. However, their proof will be completed if we replace the claim ``if $r\in R$ and $1_r.\mu = 1_\theta$, then $1_r.f(\mu) = 1_\theta$'' by ``if $r\in R$ and for some $n$, $1_{r^n}.\mu = 1_\theta$, then $1_{r^n}.f(\mu) = 1_\theta$''. We think there has been a misprint in writing. By the above arguments, we have the following theorem.
\end{rem}

\begin{thm}\label{q}\cite[Theorem 3.14]{HM2}
Let $f : M\rightarrow N$ be a monomorphism of non-zero $R$-modules such that $N/Im(f)$ is a torsion free $R$-module. If $\mu$ is a $\xi$-fuzzy second submodule of $M$, then so is $f(\mu)$.
\end{thm}

\begin{thm}\label{r}\cite[Theorem 3.15]{HM2}
Let $\lambda,\,\mu$ be two fuzzy submodules of $M$ and $1_M$ be fuzzy second with $\mu + \lambda \supseteq 1_M$. Then either $\lambda \supseteq  1_M$ or $(1_\theta : 1_M) \supseteq (1_\theta : \mu)$.
\end{thm}

\begin{cor}\label{s}\cite[Corollary 3.16]{HM2}
Let $\zeta$ be a fuzzy ideal of $R$ and let $\mu$ be a fuzzy submodule of $M$. If $1_M$ is a fuzzy second submodule of $M$ with $\mu + (1_\theta : \zeta) \supseteq 1_M$, then either $(1_\theta : \zeta) \supseteq 1_M$ or $\mu \supseteq 1_M$.
\end{cor}

Let $\mu$ be a fuzzy prime submodule of $M$ (resp. a fuzzy prime ideal of $R$) and let $\xi = (\mu : 1_M)$ (resp. $\xi = (\mu : 1_R)$). Then $\mu$ is called a \textit{$\xi$-fuzzy  prime submodule} (resp. a \textit{$\xi$-fuzzy  prime ideal}).

\begin{thm}\label{t}\cite[Theorem 3.17]{HM2}
Let $1_M \neq \lambda \in FS(M)$ and let $\zeta \in FI(R)$. Then $\zeta$ is a fuzzy prime ideal of $R$ and $\lambda$ is a $\zeta$-fuzzy prime submodule if and only if it satisfies the following conditions:
\item(a) $\forall t, s \in I, r\in R, x \in M$, if $t_r.s_x \in \lambda$, and $s_x \not\in \lambda$, then $t_r \in \zeta$.
\item(b) $\zeta \subseteq (\lambda : 1_M)$.
\end{thm}

\begin{thm}\label{u}\cite[Theorem 3.18]{HM2}
Let $1_{\theta} \neq \mu \in FS(M)$ and $\xi \in FI(R)$. Then $\xi$ is a fuzzy prime ideal of $R$ and $\mu$ is $\xi$-fuzzy second if and only if
it satisfies the following conditions:
\item(a) $r \in R$ and $1_r.\mu \neq \mu$ implies that $1_r \in \xi$.
\item(b) $\xi \subseteq (1_\theta : \mu)$.
\end{thm}

\begin{thm}\label{v}\cite[Theorem 3.19]{HM2}
Let $1_\theta \neq \mu \in FS(M)$ be $\xi$-fuzzy second, where $\xi = (1_\theta : \mu)$. Let $\nu \in FS(M)$ be such that $\nu \subset \mu$ and $\nu^* \subset \mu^*$. Then $(1_{v^*} : 1_{\mu^*})$ is a $\xi$-fuzzy prime ideal of $R$. Moreover, $(\nu : \mu)$ is a fuzzy prime ideal.
\end{thm}

\begin{thm}\label{w}\cite[Theorem 3.20]{HM2}
Let $\mu$ be a fuzzy submodule of $M$. If $(1_\theta : \mu)$ is a maximal fuzzy ideal of $R$, then $\mu$ is a  fuzzy second submodule of $M$.
\end{thm}

\begin{cor}\label{x}\cite[Corollary 3.21]{HM2}
Suppose that $\mu$ is a non-zero fuzzy submodule of $M$ which is contained in a fuzzy submodule $\lambda$ such that $(1_\theta : \lambda)$ is a maximal fuzzy ideal. Then $\mu$ is a fuzzy second submodule of $M$.
\end{cor}

\begin{cor}\label{y}\cite[Corollary 3.22]{HM2}
Let $\xi$ be a maximal fuzzy ideal such that $1_\theta \neq (1_\theta : \xi)$. Then $(1_\theta : \xi)$ is a fuzzy second submodule of $M$.
\end{cor}

\begin{defn}\cite[Definition 3.23]{HM2}
Let $\mu$ be a fuzzy submodule of $M$. Then $\mu$ is called a \emph{minimal fuzzy submodule} of $M$ if $\mu$ is a minimal element in the set of all non-zero fuzzy submodule of $M$ under pointwise partial ordering.
\end{defn}

\begin{thm}\label{z}\cite[Theorem 3.24]{HM2}
Let $\mu$ be a minimal fuzzy submodule of $M$. Then the fuzzy ideal $(1_\theta : \mu)$ is a maximal fuzzy ideal of $R$.
\end{thm}

\begin{cor}\cite[Corollary 3.25]{HM2}
Every minimal fuzzy submodule is a fuzzy second submodule.
\end{cor}
\section{$\psi$-second submodule}
\noindent
\begin{defn}\label{d2.1}\cite[Definition 2.1]{MR4401392}
Let $M$ be an $R$-module, $S(M)$ be the set of all
submodules of $M$,  and let $\psi: S(M) \rightarrow S(M) \cup \{\emptyset \}$ be a function. We say that a  non-zero submodule
$N$ of $M$ is a \textit{$\psi$-second submodule of $M$} if $r \in R$, $K$ a submodule of $M$, $rN\subseteq K$,  and $r\psi(N) \not \subseteq K $, then $N \subseteq K$ or $rN=0$.
\end{defn}

We use the following functions $\psi: S(M) \rightarrow S(M) \cup \{\emptyset \}$.
$$\psi_{M}(N)=M, \qquad \forall N \in S(M),$$
$$\psi_i(N)=(N:_MAnn_R^i(N)), \ \forall N \in S(M), \ \forall i \in \Bbb N,$$
$$\psi_\sigma(N)=\sum^{\infty}_{i=1}\psi_i(N), \qquad \forall N \in S(M).$$
Then it is clear that $\psi_M$-second submodules are weak second submodules. Clearly,  for any submodule and every positive integer $n$, we have the
following implications:
$$
second \Rightarrow \psi_{n-1}-second \Rightarrow \psi_{n}-second \Rightarrow \psi_\sigma-second.
$$
For functions $\psi, \theta: S(M) \rightarrow S(M) \cup \{\emptyset \}$, we write $\psi \leq \theta $  if $\psi(N) \subseteq \theta (N)$ for each $N \in S(M)$. So whenever $\psi \leq \theta$, any $\psi$-second submodule is $\theta$-second.

\begin{thm}\label{t2.3}\cite[Theorem 2.3]{MR4401392}
Let $M$ be an R-module and  $\psi: S(M) \rightarrow S(M) \cup \{\emptyset \}$ be a function.
Let $N$ be a $\psi$-second submodule of $M$ such that  $Ann_R(N)\psi(N)\not\subseteq N $. Then $N$ is a second submodule of $M$.
\end{thm}

\begin{cor}\label{c2.4}\cite[Corollary 2.4]{MR4401392}
Let $N$ be a weak second submodule of an $R$-module $M$ such that $Ann_R(N)M \not\subseteq N$.
Then $N$ is a second submodule of $M$.
\end{cor}

\begin{cor}\label{c42.3}\cite[Corollary 2.5]{MR4401392}
Let $M$ be an R-module and  $\psi: S(M) \rightarrow S(M) \cup \{\emptyset \}$ be a function.
If $N$ is a $\psi$-second submodule of $M$ such that  $(N:_MAnn^2_R(N))\subseteq\psi(N)$, then $N$ is a $\psi_\sigma$-second submodule of $M$.
\end{cor}

\begin{thm}\label{t2.5}\cite[Theorem 2.6]{MR4401392}
Let $M$ be an R-module and $\psi: S(M) \rightarrow S(M) \cup \{\emptyset \}$ be a function.
Let $H$ be a submodule of $M$ such that far all ideals $I$ and $J$ of $R$, $(H:_MI)\subseteq (H:_MJ)$ implies that $J \subseteq I$. If $H$ is not a second submodule of $M$, then $H$ is not a $\psi_1$-second submodule of $M$.
\end{thm}

\begin{cor}\label{c2.6}\cite[Corollary 2.7]{MR4401392}
Let $M$ be an R-module and $\psi: S(M) \rightarrow S(M) \cup \{\emptyset \}$ be a function.
Let $H$ be a submodule of $M$ such that far all ideals $I$ and $J$ of $R$, $(H:_MI)\subseteq (H:_MJ)$ implies that $J \subseteq I$. Then $H$ is a second submodule of $M$ if and only if $H$ is a
$\psi_1$-second submodule of $M$.
\end{cor}

\begin{thm}\label{t2.7}\cite[Theorem 2.8]{MR4401392}
Let $M$ be an R-module,  $\phi: S(R) \rightarrow S(R) \cup \{\emptyset \}$, and $\chi: S(M) \rightarrow S(M) \cup \{\emptyset \}$ be  functions such that $\chi(P)=\phi((P:_RM))M$.
\begin{itemize}
  \item [(a)] If $P$ is a $\chi$-prime submodule of $M$ such that $(\chi(P):_RM)\subseteq \phi((P:_RM))$, then $(P:_RM)$ is a $\phi$-prime ideal of $R$.
  \item [(b)] If $M$ is a multiplication $R$-module and $(P:_RM)$ is a $\phi$-prime ideal of $R$, then $P$ is a $\chi$-prime submodule of $M$.
\end{itemize}
\end{thm}

\begin{thm}\label{p92.8}\cite[Theorem 2.9]{MR4401392}
Let $M$ be an R-module and  $\psi: S(M) \rightarrow S(M) \cup \{\emptyset \}$, $\phi: S(R) \rightarrow S(R) \cup \{\emptyset \}$ be  functions.
\begin{itemize}
  \item [(a)] If $S$ is a $\psi$-second submodule of $M$ such that $Ann_R(\psi(S))\subseteq \phi(Ann_R(S))$, then $Ann_R(S)$ is a $\phi$-prime ideal of $R$.
 \item [(b)] If $M$ is a comultiplication $R$-module, $S$ is a submodule of $M$ such that $\psi(S)=(0:_M\phi(Ann_R(S))$,  and $Ann_R(S)$ is a $\phi$-prime ideal of $R$, then $S$ is a $\psi$-second submodule of $M$.
\end{itemize}
 \end{thm}

The following example shows that the condition ``$M$ is a comultiplication $R$-module" in   Theorem \ref{p92.8} (b) can not be omitted.
\begin{ex}\label{e22.14}\cite[Example 2.10]{MR4401392}
Let $R= \Bbb Z$, $M=\Bbb Z\oplus \Bbb Z$, and $S=2\Bbb Z \oplus 2\Bbb Z$. Clearly,  $M$ is not a comultiplication $R$-module. Suppose that  $\phi: S(R) \rightarrow S(R) \cup \{\emptyset \}$ and $\psi: S(M) \rightarrow S(M) \cup \{\emptyset \}$ be  functions such that $\phi(I)=I$ for each ideal $I$ of $R$ and  $\psi(S)=M$. Then clearly,  $Ann_R(S)=0$ is a  $\phi$-prime ideal of $R$ and  $\psi(S)=M=(0:_M\phi(Ann_R(S))$. But as
 $3S \subseteq 6\Bbb Z \oplus  6\Bbb Z$,  $S \not\subseteq 6\Bbb Z \oplus  6\Bbb Z$,  and $3S\not=0$, we have that $S$ is not a $\psi$-second submodule of $M$.
\end{ex}

\begin{prop}\label{p2.9}\cite[Proposition 2.11]{MR4401392}
Let $M$ be an R-module,  $\psi: S(M) \rightarrow S(M) \cup \{\emptyset \}$  be a function, and  $N$ be a $\psi$-second
submodule of $M$. Then we have the following statements.
\begin{itemize}
  \item [(a)] If $K$ is a submodule of $M$ with $K \subset N$ and  $\psi_K: S(M/K) \rightarrow S(M/K) \cup \{\emptyset \}$ be a function such that $ \psi_K(N/K)= \psi(N)/K$, then $N/K$ is a $ \psi_K$-second submodule of $M/K$.
  \item [(b)] Let $N$ be a finitely generated submodule of $M$, $S$ be a  multiplicatively closed subset of $R$ with $Ann_R(N) \cap S=\emptyset$, and  $S^{-1}\psi: S(S^{-1}M) \rightarrow S(S^{-1}M) \cup \{\emptyset \}$ be a function such that $(S^{-1}\psi)(S^{-1}N)= S^{-1}\psi(N)$. Then $S^{-1}N$ is a $S^{-1}\psi$-second submodule of $S^{-1}M$.
\end{itemize}
\end{prop}

\begin{prop}\label{p2.10}\cite[Proposition 2.12]{MR4401392}
Let $M$ and $\acute{M}$ be $R$-modules and  $f : M\rightarrow \acute{M}$ be an $R$-monomorphism. Let  $\psi: S(M) \rightarrow S(M) \cup \{\emptyset \}$ and $\acute{\psi}: S(\acute{M}) \rightarrow S(\acute{M}) \cup \{\emptyset \}$  be functions such that $\psi(f^{-1}(\acute{N}))=f^{-1}(\acute{\psi}(\acute {N}))$, for each submodule $\acute {N}$ of $\acute{M}$. If
 $\acute{N}$ is a $\acute{\psi}$-second submodule of $\acute{M}$ such that $\acute{N} \subseteq Im(f)$, then $f^{-1}(\acute{N})$ is a  $\psi$-second submodule of $M$.
\end{prop}

\begin{rem}\label{r22.2}\cite[Remark 2.13]{MR4401392}
Let $N$ and $K$ be two submodules of an $R$-module $M$. To prove $N\subseteq K$, it is enough to show that if $L$ is a completely irreducible submodule of $M$ such that $K\subseteq L$, then $N\subseteq L$.
\end{rem}

\begin{prop}\label{p2.11}\cite[Proposition 2.14]{MR4401392}
Let $M$ be an R-module,  $\psi: S(M) \rightarrow S(M) \cup \{\emptyset \}$  be a function, and  let $N$ be a $\psi_1$-second
submodule of $M$. Then we have the following statements.
\begin{itemize}
\item [(a)] If for $a \in R$, $aN\not=N$, then $(N:_MAnn_R(N))\subseteq (N:_Ma)$.
\item [(b)] If $J$ is an ideal of $R$ such that $Ann_R(N)\subseteq J$ and $JN \not=N$, then
$(N:_MAnn_R(N))=(N:_MJ)$.
\end{itemize}
\end{prop}

\begin{thm}\label{t2.12}\cite[Theorem 2.15]{MR4401392}
Let $M$ be an $R$-module,  $\psi: S(M) \rightarrow S(M) \cup \{\emptyset \}$  be a function, and let $a$ be an element of $R$ such that $(0 :_M a) \subseteq  a(0:_MaAnn_R((0:_Ma)))$. If $(0:_Ma)$ is a $\psi_1$-second submodule of $M$, then  $(0:_Ma)$ is a second submodule of $M$.
\end{thm}

\begin{thm}\label{t2.13}\cite[Theorem 2.16]{MR4401392}
Let $N$ be a non-zero submodule of an $R$-module $M$ and  $\psi: S(M) \rightarrow S(M) \cup \{\emptyset \}$  be a function.
Then the following are equivalent:
\begin{itemize}
\item [(a)] $N$ is a $\psi$-second submodule of $M$;
\item [(b)] for completely irreducible submodule $L$ of $M$ with $N \not \subseteq L$, we have $(L:_RN)=Ann_R(N) \cup (L:_R\psi(N))$;
\item [(c)] for completely irreducible submodule $L$ of $M$ with $N \not \subseteq L$, we have $(L:_RN)=Ann_R(N)$ or $(L:_RN)=(L:_R\psi(N))$;
\item [(d)] for any ideal $I$ of $R$ and any submodule $K$ of $M$, if $IN \subseteq K$ and $I\psi(N) \not \subseteq K$,
then $IN=0$ or $N \subseteq K$.
\item [(e)] for each $a \in R$ with $a\psi(N) \not \subseteq aN$, we have $aN=N$ or $aN=0$.
\end{itemize}
\end{thm}

\begin{ex}\label{e2.114}\cite[Example 2.17]{MR4401392}
Let $N$ be a non-zero submodule of an $R$-module $M$ and  let $\psi: S(M) \rightarrow S(M) \cup \{\emptyset \}$  be a function.
If $\psi(N)=N$, then  $N$ is a $\psi$-second submodule of $M$ by Theorem \ref{t2.13} $(e)\Rightarrow (a)$.
\end{ex}

Let $R_1$ and $R_2$ be two commutative rings with identity. Let $M_1$ and $M_2$ be $R_1$ and
$R_2$-module, respectively and put $R = R_1 \times R_2$. Then $M = M_1 \times M_2$ is an $R$-module
and each submodule of $M$ is of the form $N = N_1 \times N_2$ for some submodules $N_1$ of
$M_1$ and $N_2$ of $M_2$. Suppose that $\psi^i: S(M_i) \rightarrow S(M_i) \cup \{\emptyset \}$  be a function for $i=1, 2$. The second submodules of the $R = R_1 \times R_2$-module $M = M_1 \times M_2$ are in the form $S_1 \times 0$ or $0\times S_2$, where $S_1$ is a second
submodule of $M_1$ and $S_2$ is a second submodule of $M_2$ \cite[Lemma 2.23]{AF1112}.  The following example, shows that this is not true for correspondence
$\psi^1 \times \psi^2$-second submodules in general.

\begin{ex}\label{e2.14}\cite[Example 2.18]{MR4401392}
Let $R_1 = R_2 = M_1 = M_2 =S_1= \Bbb Z_6$. Then clearly,  $S_1$  is
a weak second submodule of $M_1$. However,
$
(\bar{2},\bar{1})(\Bbb Z_6 \times 0) \subseteq \bar{2}\Bbb Z_6 \times \bar{3}\Bbb Z_6
$
and $(\bar{2},\bar{1})(\Bbb Z_6 \times \Bbb Z_6) \not \subseteq \bar{2}\Bbb Z_6 \times \bar{3}\Bbb Z_6$. But
$(\bar{2}, \bar{1})(\Bbb Z_6 \times 0) =\bar{2}\Bbb Z_6\times 0 \not =0\times 0$, and $\Bbb Z_6 \times 0 \not \subseteq \bar{2}\Bbb Z_6 \times \bar{3}\Bbb Z_6$.
Therefore,  $S_1 \times 0$ is not a weak second submodule of $ M_1 \times M_2$.
\end{ex}

\begin{thm}\label{t2.133}\cite[Theorem 2.19]{MR4401392}
Let $R = R_1 \times R_2$ be a ring and $M = M_1 \times M_2$
be an $R$-module, where $M_1$ is an $R_1$-module and $M_2$ is an $R_2$-module. Suppose that $\psi^i: S(M_i) \rightarrow S(M_i) \cup \{\emptyset \}$  be a function for $i=1, 2$. Then  $S_1 \times 0$ is a $\psi^1 \times \psi^2$-second submodule of $M$,  where $S_1$ is a  $\psi^1$-second submodule of $M_1$ and  $\psi^2(0)=0$.
\end{thm}

Let $\psi_M : S(M)\rightarrow S(M)$ be the function defined by $\psi_M(L) = M$ for every $L \in S(M)$. Then a $\psi_M$-second
submodule of $M$ is said to be a \textit{weak second submodule} of $M$.
Let $n \geq 2$ be an integer and $\psi_n : S(M)\rightarrow S(M)$ be the function defined by $\psi_n(L) = (L :_M Ann_R(L)^{n-1})$ for every $L \in S(M)$. Then a $\psi_n$-second submodule of $M$ is said to be an \textit{$n$-almost second submodule} of $M$. In particular, for $n = 2$, a $2$-almost second submodule of $M$ is called an \textit{almost second submodule} of $M$.
Let $M$ be an $R$-module and $N$ be a submodule of $M$. Since $\psi(N)\setminus N = (\psi(N) \cup N)\setminus N$, without loss of
generality, throughout this paper we will assume that $N \subseteq \psi(N)$.
It is clear from the definition that every $R$-module $M$ is a $\psi$-second submodule of itself for any function
$\psi : S(M)\rightarrow S(M)$. But not every $R$-module is a second submodule of itself. For example, $\Bbb Z$ is not a second
$\Bbb Z$-submodule of itself but $\Bbb Z$ is $\psi $-second $\Bbb Z$-submodule of itself.

\begin{thm}\label{tuu2.134}\cite[Theorem 2.3]{MR4290828}
Let $N$ be a $\psi$-second submodule of an $R$-module  $M$ such that $N = M_1 + \cdots + M_k$, where $M_1,\cdots ,M_k $ are submodules
of $M$ such that $Ann_R(M_i)$ is a maximal ideal of $R$ for each $i $ ($1 \leq i \leq k$). Then either $N$ is second or $N =\psi (N)$.
\end{thm}

\begin{cor}\label{tuu2.134}\cite[Corollary 2.4]{MR4290828}
Let $M$ be an $R$-module and $N$ be a $\psi$-second submodule of $M$. If one of the following holds, then either
$N$ is second or $N =\psi (N)$.
\begin{itemize}
\item [(a)] $M$ is a finitely generated semisimple $R$-module.
\item [(b)] $M$ is a finitely cogenerated co-semisimple $R$-module.
\end{itemize}
\end{cor}

\begin{prop}\label{tuu2.135}\cite[Proposition 2.6]{MR4290828}
Let $M$ be a comultiplication $R$-module and $N$ be a non-zero submodule of $M$. Then the following
are equivalent.
\begin{itemize}
\item [(a)] $N$ is a $\psi$-second submodule of $M$.
\item [(b)] If $K$ and $L$ are two submodules of $M$ such that $N \subseteq C(KL)$ and $\psi (N)\not \subseteq C(KL)$, then $N \subseteq K$ or $N \subseteq L$.
\end{itemize}
\end{prop}

\begin{prop}\label{tuu2.136}\cite[Proposition 2.7]{MR4290828}
Let $M$ be an $R$-module and $K$ be a submodule of $M$. Let $\psi^K : S(K)\rightarrow S(K)$ be the function defined
by $\psi^K(L) = \psi (L) \cap K$ for every $L \in S(K)$. Then the following hold for a submodule $N$ of $K$.
\begin{itemize}
\item [(a)] If $N$ is a $\psi$-second submodule of $M$, then $N$ is a $\psi^K$-second submodule of $K$.
\item [(b)] Let $\psi (N) \subseteq K$. Then $N$ is a $\psi^K$-second submodule of $K$ if and only if $N$ is a $\psi$-second submodule of $M$.
\item [(c)] If $K \subseteq \psi (N)$ and $N$ is a $\psi$-second submodule of $M$, then $N$ is a weak second submodule of $K$.
\item [(d)] If $\psi (N) \subseteq \psi (K)$, $K$ is a $\psi$-second submodule of $M$ and $N$ is a weak second submodule of $K$, then $N$ is a $\psi$-second submodule of $M$.
\end{itemize}
\end{prop}

\begin{prop}\label{tuu2.138}\cite[Proposition 2.8]{MR4290828}
Let $R_i$ be a ring and $M_i$ be an $R_i$-module for $i = 1, 2$. Denote $R = R_1 \times R_2$ and $M = M1 \times M2$.
Suppose that $N_1$ is a weak second submodule of $M_1$ such that $\psi (N1 \times \{0\}) \subseteq M1 \times \{0\}$. Then $N1 \times \{0\}$ is a $\psi$-second
submodule of $M$.
\end{prop}

\begin{cor}\label{tuu2.137}\cite[Corollary 2.9]{MR4290828}
Let $R_i$ be a ring, $M_i$ be an $R_i$-module for $i = 1, 2$. Denote $R = R_1 \times R_2$ and $M = M_1 \times M_2$. Suppose
that $\psi : S(M) \rightarrow S(M)$ is a function with $\psi \leq  \psi_\omega$. Then, $N_1 \times \{0\}$ is a $\psi$-second submodule of $M$ for any weak second submodule of $N_1$ of $M_1$.
\end{cor}

\begin{prop}\label{tuu2.139}\cite[Proposition 3.1]{MR4290828}
Let $M$ be a prime $R$-module and $N$ be a proper almost second submodule of $M$. Then $N$ is a prime
submodule of $M$.
\end{prop}

\begin{thm}\label{tuu2.1310}\cite[Theorem 3.2]{MR4290828}
Let $M$ be an $R$-module. If there exist maximal ideals $P1, \ldots,  Pn$ of $R$ such that $P_1 \cap \cdots \cap P_n \subseteq Ann_R(M)$, then for any ideal $I$ of $R$, $(0 :_M I) = 0$ or $(0 :_M I)$ is an almost second submodule of $M$.
\end{thm}

\begin{thm}\label{tuu2.1311}\cite[Theorem 3.3]{MR4290828}
Let $M$ be an $R$-module, $a \in R$, $(0 :_M a) \not= (0)$ and $(0 :_M a) = aM$. Then $(0 :_M a)$ is an almost second
submodule of $M$ if and only if it is a second submodule of $M$.
\end{thm}

\begin{lem}\label{tuu2.1312}\cite[Lemma 3.4]{MR4290828}
An $R$-module $M$ is fully coidempotent if and only if $N= (N :_M Ann_R(N)^m)$ for every submodule $N$ of
$M$ and positive integer $m$.
\end{lem}

\begin{thm}\label{tuu2.1313}\cite[Theorem 3.5]{MR4290828}
 Let $R = R_1 \times \cdots \times R_m$ and $M = M_1 \times \cdots \times M_m$ where $R_i$ is a ring, $0\not= M_i$ is an $R_i$-module for all $i \in \{1, \ldots ,m \}$ and $n,m \geq 2$. Then every non-zero submodule of $M$ is $n$-almost second if and only if $M$ is a fully
coidempotent $R$-module.
\end{thm}

\begin{cor}\label{tuu2.1314}\cite[Corollary 3.6]{MR4290828}
 Let $R = R_1 \times \cdots \times R_m$ and $M = M_1 \times \cdots \times M_m$ where $R_i$ is a ring, $0\not= M_i$ is an $R_i$-module for all $i \in \{1, \ldots ,m \}$ and $n,m \geq 2$.  Then every non-zero submodule of $M$ is $n$-almost second if and only if every non-zero
submodule of $M$ $(n+1)$-almost second.
\end{cor}

\begin{lem}\label{tuu2.1315}\cite[Lemma 3.7]{MR4290828}
Let $M$ be an $R$-module, $N$ be a submodule of $M$ and $I$ be an ideal of $R$. Suppose that $(0 :_M I) \not= (N :_M I)$
and $(N :_M I) \not= N$. Then $K := (N :_M I)$ is an almost second submodule of $M$ if and only if $K = (K :_M Ann_R(K))$.
\end{lem}

\begin{thm}\label{tuu2.1316}\cite[Theorem 3.8]{MR4290828}
Let $M$ be an Artinian $R$-module, $I \subseteq Jac(R)$ and $N$ be a submodule of $M$ such that $(N :_R M) = 0$ and
$(0 :_M I) \not=(N :_M I)$. Then $(N :_M I)$ is not an $n$-almost second submodule of $M$ for any integer $n > 1$.
\end{thm}

\begin{lem}\label{tuu2.1317}\cite[Lemma 3.9]{MR4290828}
Let $I$ be an ideal of $R$ and $M$ be an $R$-module. Then $(0 :_M (Ann_R(0 :_M I))^n) = (0 :_M I^n)$ for every integer
$n > 1$. In particular, $((0 :_M I) :_M (Ann_R(0 :_M I))^{n-1}) = (0 :_M I^n)$.
\end{lem}

\begin{prop}\label{tuu2.1318}\cite[Proposition 3.10]{MR4290828}
Let $M$ be a strong comultiplication $R$-module and $I$ be an ideal of $R$. Then, $I$ is an $n$-almost prime
ideal of $R$ if and only if $(0 :_M I)$ is an $n$-almost second submodule of $M$.
\end{prop}

\begin{ex}\label{tuu2.1319}\cite[Example 3.11]{MR4290828}
Let $R = F[[X^3,  X^4, X^5]]$ where $F$ is a field, and $I = RX^3 + RX^4$. Then $I$ is an almost prime ideal of $R$
which is not a $3$-almost prime ideal by \cite[Example 11]{MR2388034}. Let $M$ be a strong comultiplication $R$-module. By Proposition \ref{tuu2.1318}, $(0 :_M I)$ is an almost second submodule of $M$ which is not a $3$-almost second submodule of $M$.
\end{ex}

\begin{ex}\label{tuu2.1321}\cite[Example 3.12]{MR4290828}
Let $M$ be a co-semisimple $R$-module. Then \cite[Theorem 4.8]{AF101} and  \cite[Theorem 2.3]{MR2101067} imply that
$(0 :_M I) = (0 :_M I^n)$ for every ideal $I$ of $R$ and an integer $n > 1$. By Lemma \ref{tuu2.1317}, $(0 :_M I)$ is an $n$-almost second submodule of $M$ for each integer $n > 1$.
\end{ex}

\begin{thm}\label{tuu2.1321}\cite[Theorem 3.13]{MR4290828}
Let $n > 1$ be an integer, $M$ be an $R$-module and $I$ be an ideal of $R$ with $(0 :_M I) \not=0$.
\begin{itemize}
\item [(a)] If $R$ is a $ZPI$-ring and $(0 :_M I)$ is an $n$-almost second submodule of $M$, then $(0 :_M I) = (0 :_M I^n)$ or
$(0 :_M I) = (0 :_M P)$ for some prime ideal $P$ of $R$.
\item [(b)] If $R$ is a Dedekind domain, then $(0 :_M I$) is an $n$-almost second submodule of $M$ if and only if $(0 :_M I) = (0 :_M I^n)$ or $(0 :_M I)$ is a second submodule of $M$.
\item [(c)] If $(R, \mathfrak{m})$ is a local $ZPI$-ring and $(0 :_M I)$ is finitely cogenerated, then $(0 :_M I)$ is an $n$-almost second submodule of $M$ if and only if $(0 :_M I) = M$ or $(0 :_M I) = (0 :_M\mathfrak{m})$.
\end{itemize}
\end{thm}

\begin{thm}\label{tuu2.1322}\cite[Theorem 3.14]{MR4290828}
Let $N$ be an $n$-almost second submodule of an $R$-module $M$ and $I = Ann_R(N)$. Then $S = [(R\setminus I) \cup
Ann_R(N :_M I^{n-1})]\setminus P$ is a $P$-essential multiplicatively closed subset of $R$ for each prime ideal $P$ of $R$.
\end{thm}

\begin{defn}\label{tuu2.1323}\cite[Definition 3.15]{MR4290828}
Let $M$ be an $R$-module, $N$ be a submodule of $M$ and
$T := \{Q \leq N : Q\  is\ almost\ second\ and \  (N :_M Ann_R(N)) = (Q :_M Ann_R(Q))\}$.
Then almost second radical of $N$ is defined as the submodule $a$-$sec(N) :=\sum_{Q\in T} Q$ if $T\not=0$. If $T =0$ then $a$-$sec(N)$ is defined as $(0)$.
\end{defn}

\begin{defn}\label{tuu2.1324}\cite[Definition 3.16]{MR4290828}
Let $M$ be an $R$-module and $S$ be a proper subset of $M$. If, for any submodules $K$, $L$ of $M$ and any ideal
$I$ of $R$, $(K \cap L)\cup S \not =M$, $(K \cap (0 :_M I))\cup S\not=M$ and $(S_c :_M Ann_R(S^c))\not \subseteq (L :_M I)$ imply that $(K \cap (L :_M I)) \cup S \not= M$, then $S$ is called an \textit{almost $m^*$-system}.
\end{defn}

\begin{prop}\label{tuu2.1324}\cite[Proposition 3.17]{MR4290828}
Let $M$ be an $R$-module, $Q$ be a non-zero submodule of $M$. Then, $Q$ is an almost second submodule
of $M$ if and only if $S := M\setminus Q$ is an almost $m^*$-system.
\end{prop}

\begin{thm}\label{tuu2.1325}\cite[Theorem 3.18]{MR4290828}
Let $S$ be an almost $m^8$-system in $M$ and $Q$ be a non-zero submodule of $M$ minimal with respect to the
properties that $Q \cup S = M$ and $(Q :_M Ann_R(Q)) = (S^c :_M Ann_R(S^c))$. Then $Q$ is an almost second submodule of $M$.
\end{thm}

\begin{thm}\label{tuu2.1326}\cite[Theorem 3.19]{MR4290828}
 Let $M$ be an $R$-module and $N$ be a submodule of M. If there exists an almost second submodule $Q$ of
$N$ with $(N :_M Ann_R(N)) = (Q :_M Ann_R(Q))$, then $a$-$sec(N) = \{x \in N : x \not \in  S, N \cup S = M\ and\ (N :_M Ann_R(N)) = (S^c :_M Ann_R(S^c)\ for\ some\ almost\ m^*-system\ S\ in\ M\}$.
\end{thm}

\begin{Conclusion}\label{c8.31}
As we mentioned in the introduction, there is a large body of researches related to second submodules of modules since this notion has been introduced. Also, there is large open space for this notion parallel to researches on prime submodules.  In \cite{MR3648183, MR3954214},  the authors applied the notion of second submodules in lattice theory. Also, the concept of second submodules of module can be applied in other fields such as  graph theory.
\end{Conclusion}

\textbf{Conflicts of Interest.}
The author declares that there are no conflicts of interest.


\bibliographystyle{amsplain}

\end{document}